\documentclass[11pt, a4paper]{article}
\usepackage{verbatim}
\usepackage{amsmath,latexsym}
\usepackage{amssymb,wasysym}
\usepackage[latin1]{inputenc}
\usepackage{soul}%define \st para tachare
\usepackage{slashbox}
\usepackage{ marvosym }
\usepackage{array}
\usepackage{tabularx}
\usepackage{fancyhdr}
\usepackage{psfrag}
\usepackage{xspace}
\usepackage{pstcol}
\usepackage{pstricks}
\usepackage{pst-node,xcolor,graphicx}
\usepackage{pst-tree}
\usepackage{graphicx}
\usepackage{subfigure}
\usepackage{color}
\usepackage{colortbl}
\usepackage{rotating}
\usepackage{multirow}
\usepackage{longtable}
\usepackage{capt-of}
\usepackage[font=small]{caption}
\usepackage[round]{natbib}
\usepackage{apalike}
\usepackage{url}

\newcolumntype{Y}{>{\centering\arraybackslash}X}

\definecolor{maroon}{cmyk}{0,0.87,0.68,0.32}
\definecolor{lightblue}{rgb}{.68,.85,.9}

\newcommand{\dem}{\par \noindent{\bf Proof:} \par}
\newcommand{\fin}{\hfill $\square$  \par \bigskip}

\newtheorem{teor}{\bf Theorem}[section]

\newtheorem{lema}{\bf Lemma}[section]

\newtheorem{ejem}{\bf Example}[section]

\newcommand{\ABC}{\ensuremath{PpCP}\xspace}
\newcommand{\kABC}{\ensuremath{K\textrm{-}PpCP}\xspace}
\newcommand{\UB}{\ensuremath{U\!B}}
\newcommand{\LB}{\ensuremath{L\!B}}

\allowdisplaybreaks[1]

\definecolor{mari}{rgb}{0.24,0.29,0.98}

\begin{document}
\title{The probabilistic $p$-center problem:\\ Planning service for potential customers}
\author{Luisa I. Mart\'{\i}nez-Merino$^1$, Maria Albareda-Sambola$^{2}$, Antonio M.\ Rodr\'{\i}guez-Ch\'{\i}a$^1$}

\date{\scriptsize
$^1$Departamento de Estad\'{\i}stica e Investigaci\'on Operativa, Universidad de C\'adiz, Spain\\
$^2$ Departamento de Estad\'{\i}stica e
Investigaci\'on Operativa, Universitat Polit\`{e}cnica de Catalunya-BarcelonaTech, Barcelona, Spain \\{\color{white}.}}

\maketitle

\begin{abstract}
This work deals with the probabilistic p-center problem, which aims at minimizing
the expected maximum distance between any site with demand and its center, considering that each site has demand with a specific probability. The problem is of interest when emergencies may occur at predefined sites with known probabilities. For this problem we propose and analyze different formulations as well as a Variable Neighborhood Search heuristic. Computational tests are reported, showing the potentials and limits of each formulation, the impact of their enhancements, and the effectiveness of the heuristic.
\end{abstract}

\noindent{\bf Keywords:} Discrete location, P-center, Mixed integer linear formulations, {Demand uncertainty.}

\section{Introduction}
\label{sec_int}
Many discrete location models have been inspired by a variety of applications in logistics, telecommunications, emergency services, etc.  The goal is to locate a number of facilities within a set of candidate sites and assign customers to them optimizing some effectiveness measure, usually depending on the assignment distances \citep[see][]{LocationBook, daskin1995, Drezner.ea.location}.

Among them,  the $p$-center problem ($pCP$) aims at locating $p$ centers out of $n$ sites and assigning the remaining sites to the centers, so that the maximum distance between a site and its assigned center is minimized \citep[see][]{pCenterChapter}.
Although the $pCP$ is NP-hard \citep[][]{karivhakimisjam1979_1}, it can be solved  efficiently via bisection search \citep[see][]{daskin1995,daskin2000}. Nonetheless, extensive literature exists proposing exact and heuristic algorithms for \citep[][]{doublebound_pc,VNS_hybrid}.
The main applications of the $pCP$ are the location of emergency services like ambulances, hospitals or fire stations, since, in this context, the whole population should be timely reachable from some center. However, as already observed in the past, locating services according to the $pCP$ may increase the effective service distances {\citep[see][]{Ogryczak97}}. This motivated alternative models, such as the cent-dian \citep[][]{Halpern78}.

This work presents a stochastic $pCP$ variant {\citep[see, for  instance,][]{snyderdaskintr2005, bermandrassmenezesor2007,Sergio,snyderiiet2006}. This variant} aims at smoothing this loss of spatial efficiency, trying to keep the centers close to where they are needed.
Namely, the probabilistic $p$-center problem (\ABC) aims at finding $p$ centers, out of $n$ sites, that minimize the expected maximum distance between a site with demand and its allocated center, assuming that demands can occur at each site independently, and with a known probability.

As stated above, considering the expected maximum service cost instead of the maximum assignment distance, prevents situations where a remote site with a  low demand probability forces to place centers further from the remaining sites than it is desirable.
In applications like firefighting, for instance, one pretends to provide service to a whole region but it wouldn't make sense to use a worst-case approach if the region contains areas with high risk of fire, and others where a fire is very unlikely to take place.
{Another real situation that fits to our model could be the case of locating several security offices to attend as fast as possible burglary alarms in different neighborhoods of a city. Obviously, by beforehand studies and previous experience, it would be possible 
to know that some neighborhoods are safer than others, and then an estimation of the probability that a theft takes place in a
neighborhood can be  computed. In this case, the probabilities that a theft occurs in two different neighborhoods are independent. 
In such  situations, the \ABC\ would be much more convenient than the classical $pCP$.}

From the modeling point of view, the \ABC\ falls into the stochastic programming paradigm, where uncertain values are described through probability distributions \citep[see, for instance,][]{AlbaredaSambola2011335,Huang10}
 as opposite to the robust optimization  approach, which attempts to optimize the worst-case system performance when uncertain data is only described using data ranges \citep[e.g.,][]{koyu97,PuertoChia, EMR15,Lu13,LuSheu13}. The \ABC also differs from other analyzed location problems where the centers are not restricted to be nodes of a network \citep[see ][]{StocALRChapter}.

For the \ABC we explore three formulations and a variable neighborhood search (VNS) heuristic.
Within the formulations, we have considered an ordered objective function \citep[see ][]{nickelpuerto2005}.
This function weights the assignment costs with different factors that depend on their position in the ordered list of incurred costs. In the \ABC, these factors are decision variables, since each one depends on the customers that have larger costs.

The paper is organized as follows. Section \ref{sec_the_prob} defines and analyzes the \ABC and the more general \kABC, where only the $K$ largest assignment distances are considered. Section \ref{section_homogeneous} focuses on the homogeneous case (all customers share the same demand probability). The alternative formulations for the general \kABC\ and their enhancements are exposed in Section~\ref{sec_formulations}.
Lower and upper bounds are discussed in Section~\ref{bounds1} and a VNS heuristic is presented in Section~\ref{sec_vns}. The computational experiments evaluating the formulations and their enhancements, the quality of the bounds, and the efficacy of the heuristic, are reported in Section~\ref{sec_comput}. Our findings and future research lines conclude this work.

\section{The problem}
\label{sec_the_prob}
Let $N=\{1,\ldots, n\}$ be the given set of customer sites. Throughout the paper we assume, without loss of generality, that the set of candidate sites for centers is identical to $N$, although all results apply in the case where only some of them are eligible. Let $p\geqslant2$ be the number of centers to be located. For each pair $i,j\in N$,
let $d_{ij}$ be the distance (service cost) from $i$ to $j$. We assume $d_{ii}=0$ $\forall i\in N$ and $d_{ij}>0$ $\forall i\neq j \in N$ (these distances need not to be proper distances, since triangle inequality is not assumed to hold).
In case of ties among several distances from the same site we assume without loss of generality, that preferences are given by the site index. Accordingly, in what follows, site $i$ will prefer center $j$ rather than $j'$, denoted by $d_{ij}\prec d_{ij'}$, whenever $d_{ij}<d_{ij'}$ or $d_{ij}=d_{ij'}$ and $j<j'$. Finally, service requests at the customer sites are assumed to take place independently with probabilities $0< q_i \leqslant 1$, $i\in N$.

A solution to the $PpCP$ consists of a set of $p$ centers, plus the assignment of each site to one of them. However, at the moment of making the decision, we do not know which customers will indeed place a request. Therefore, once demands are revealed, only the service of customers with demand will incur a cost. Accordingly, in what follows, we will distinguish between assignment distances (distances between customers and their respective assigned centers) and service costs (distances between customers where demand occurs and their respective assigned centers). The goal of the \ABC is to identify the solution with the smallest expected value (among all scenarios) of the maximum service cost.

\begin{ejem}
Given the set of sites $N$ with coordinates $N=\{(21,39),(37,16),(19,26),(71,26),(25,59),$ $(85,39),(88,59),(82,59),(15,86),$
$(41,26)\}$, and using Euclidean distances; consider the three instances of the $P3CP$ defined by the following three probability vectors:
\begin{itemize}
\addtolength{\itemsep}{-12pt}
\item $q_1=(\:0.06\:,0.05\:,0.07\:,0.02\:,0.1\:,0.11\:,0.18\:,0.09\:,0.01\:,0.16\:)\:,$
\item $q_2=(\:0.45\:,0.56\:,0.51\:,0.46\:,0.41\:,0.54\:,0.59\:,0.43\:,0.44\:,0.52\:)$ and
\item $q_3=(\:0.89\:,0.84\:,0.82\:,0.81\:,0.83\:,0.88\:,0.83\:,0.96\:,0.94\:,0.92\:).$
\end{itemize}
The instances and the corresponding optimal solutions are shown in Figure~\ref{fig:ejemplo}. Each circle represents a site, and its size is proportional to its corresponding $q$ value. Optimal centers are filled in black.
\begin{figure}[!h]
\fbox{\includegraphics[width=0.31\textwidth]{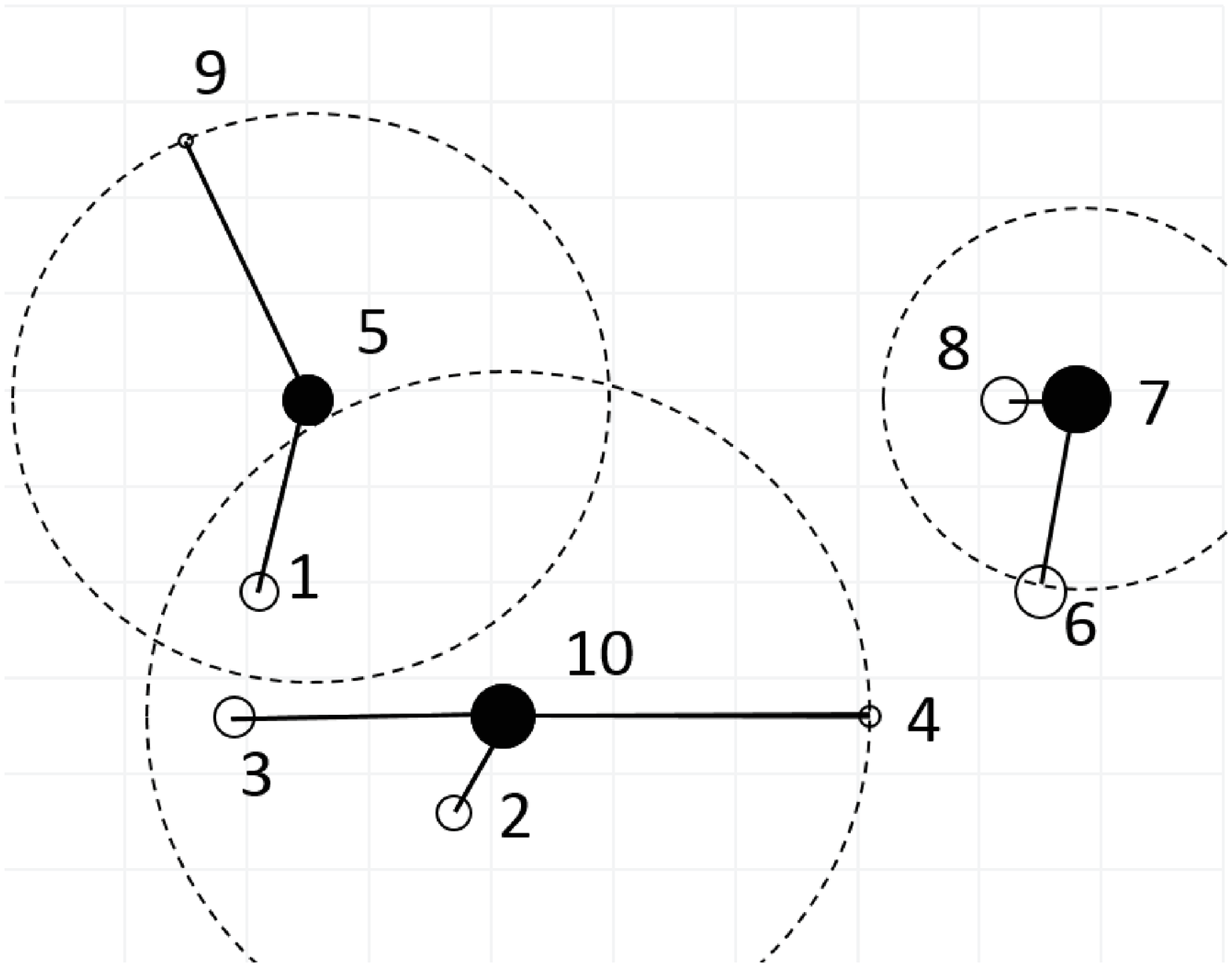}} \ \fbox{\includegraphics[width=0.31\textwidth]{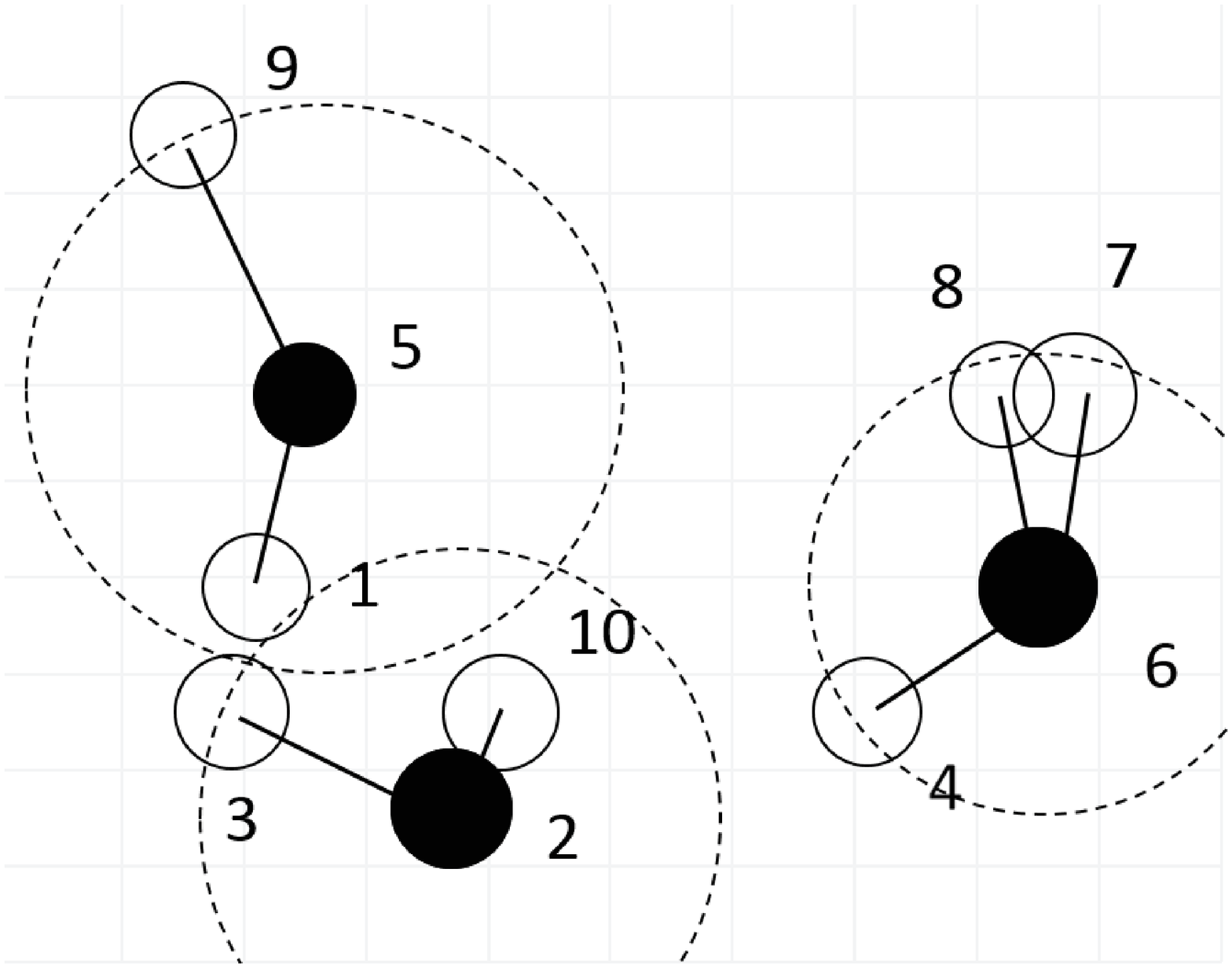}}
\ \fbox{\includegraphics[width=0.31\textwidth]{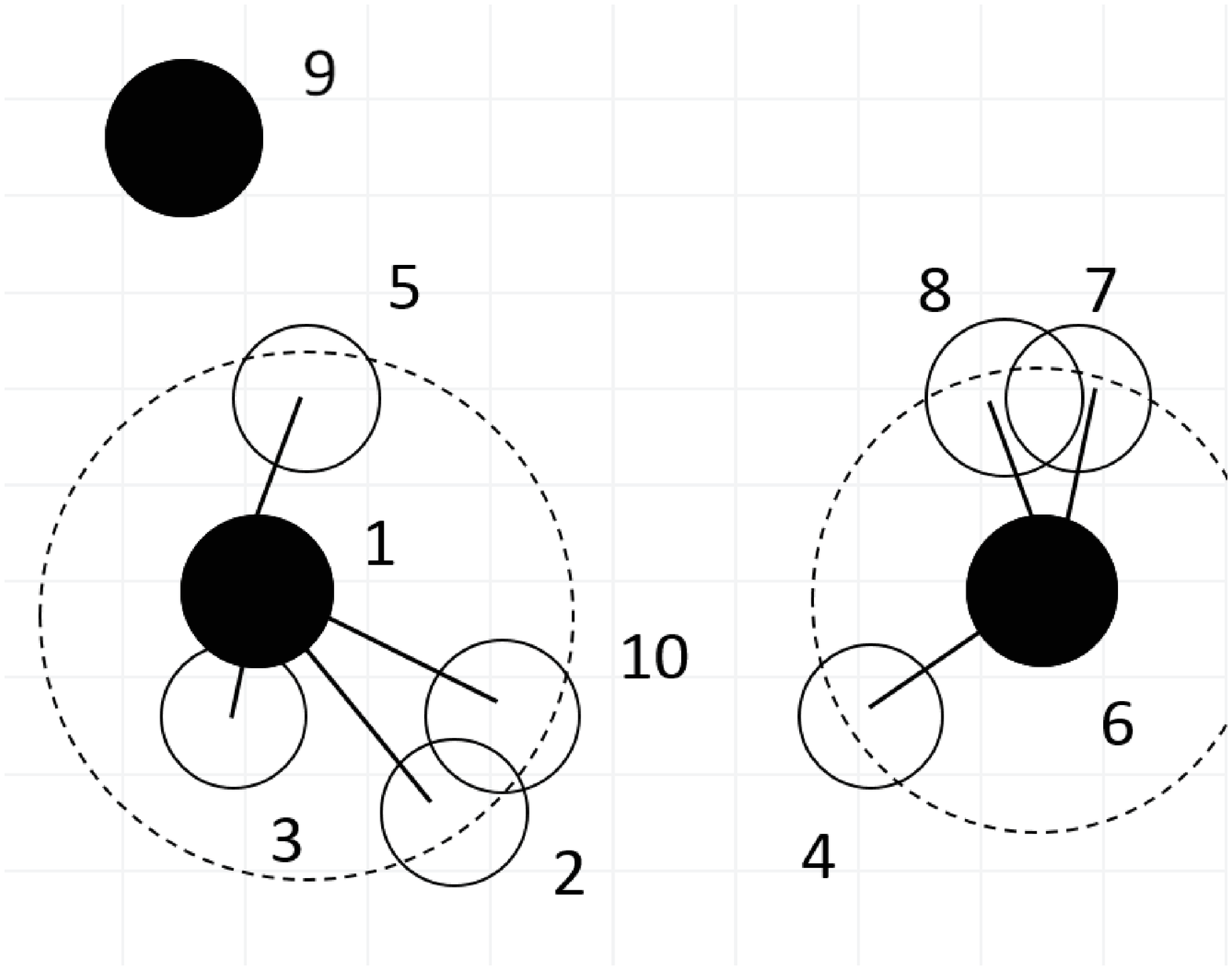}}
\centerline{\phantom{es} ${\mathbf q_1:}$ $z^*=15.69$. $\equiv 3$-median\phantom{espa} ${\mathbf q_2:}$ $z^*=23.68$. $\equiv 4$-$3$-centrum \phantom{espa} $\mathbf{q_3}$: $z^*=27.31$. $\equiv 3$-center\phantom{es}}
\caption{ Solutions with different demand probabilities\label{fig:ejemplo}}
\end{figure}
As can be observed, when demand probabilities are small ($q_1$), the optimal centers for the $P3CP$ coincide with the optimal solution of the $3$-median problem. Similarly, the solution of the $4$-$3$-centrum  (locating $3$ facilities with the $4$-centrum criterion)
is optimal for the $P3CP$ with demand probabilities given by $q_2$.
Finally, the $P3CP$ and the $3CP$ have the same solution for large demand probabilities ($q_3$).
\end{ejem}

The above example illustrates the typical behavior of the \ABC in relation to classical location models, for different $q$ values. Indeed, if demand probabilities are similar and very small, the probabilities of each assignment distance yielding the largest service cost become very similar and, therefore, the \ABC resembles the $p$-median problem. As opposite, if these probabilities are high, the probability that the furthest assignment yields the largest service cost is almost 1 and, therefore, all the other assignment distances have small weights in the objective function, leading to solutions similar to those of the $pCP$. That is, depending on the demand probabilities, the \ABC may yield a whole range of solutions. Therefore, the \ABC can be seen as a tradeoff between classical discrete location models that focus on reducing the largest assignment distances, such as the $pCP$ or the $k$-centrum, and those that minimize the total service distance, like the $p$-median. Analogously, from the managerial point of view, the model presented here allows to identify solutions that represent a tradeoff between the quality of service (associated with minimizing the largest assignment distance) and the cost of service (associated with minimizing the total assignment cost).

Recall that the objective function of the \ABC  accounts for the expected maximal service cost.
To compute this expected value for a solution where the set of located centers is $J\subset N$, we will use a matrix
$(\pi_{ij})_{i\in N,j\in N}$. Hence, if site $i$ is assigned to a center located at $j\in J$,
$\pi_{ij}$ will be the probability that there is no demand at the sites whose assignment distances are larger than $d_{ij}$, and it will take value $0$ otherwise.

\begin{lema} \label{lemapis}
For a given solution with  centers located at $J\subset N$, the matrix   $(\pi_{ij})_{i\in N,j\in N}$ satisfies:
\begin{enumerate}
\item $| \{j\in J: \pi_{ij}\neq 0 \}| \leqslant 1$  and $\pi_{ij}=0\;\forall j\notin J,\forall i\in N$.
\item Let $d_{(1)} \leqslant \cdots \leqslant d_{(n)}$ be a non-decreasing sequence of assignment distances and $(1), \ldots, (n)$ the corresponding sequence of customers. For $i\leqslant n$,
     $\sum_{j\in N}\pi_{(i)j}=\sum\nolimits_{ j\in J} \pi_{(i)j} = \prod_{t=i+1}^{n} (1 - q_{(t)}) .$
\item The expected maximum service cost can be computed as  $\sum_{i=1}^n \sum_{j=1}^n \pi_{ij} q_i  d_{ij}=\sum_{i=1}^n \sum_{ j\in J} \pi_{ij} q_i  d_{ij}.$
\item It holds that
\begin{equation}
\sum_{i=1}^n \sum_{j\in J} q_i\pi_{ij} = 1- \prod\limits_{j=1}^{n} (1-q_{j})
\leqslant 1.\label{eq:sumapis}\end{equation}
\end{enumerate}
\end{lema}

\dem

\begin{enumerate}
\addtolength{\itemsep}{-12pt}
\item Follows from the single assignment assumption and the definition of $\pi$.
\item Given a solution, for each $i\in N$ let $j_{i}$ be its assigned center in the solution. Then, by 1), $ \sum\nolimits_{ j\in J} \pi_{(i)j} = \pi_{(i)j_{(i)}}$. Now, by definition, $\pi_{(i)j_{(i)}}
    = \prod_{t=i+1}^{n} (1 - q_{(t)})$; that is, the probability that all sites with assignment costs larger than $d_{(i)j_{(i)}}$ have no demand. Note that this can be computed as the product for all these sites of the probability of not having demand, since service requests are assumed to be independent.
\item Note that, a given assignment distance $d_{ij_i}$ will become a service cost only if \textsl{i)} site $i$ has demand (which happens with probability $q_i$); and \textsl{ii)} no site with a larger assignment distance does (which happens with probability $\pi_{ij_i}$). Therefore, the expected service cost can be computed as
    $\sum\nolimits_{i=1}^n (q_i\pi_{ij_i}) d_{ij_i}$.
    Since $\pi_{ij}=0\ \forall j\neq j_i$, all the other terms in \textsl{3)} are zero and the result holds.
\item  $\sum\nolimits_{i=1}^n \sum_{j\in J} q_i\pi_{ij} =\sum\nolimits_{i=1}^n q_i\pi_{ij_i}= q_{(1)}+q_{(2)}\pi_{(2)j_{(2)}}+ \cdots q_{n}\pi_{(n)j_{(n)}}$. This is exactly the probability that at least one site has demand. The complement of this event consists of the single scenario where no site has demand, which has probability $\prod\nolimits_{j=1}^{n} (1-q_{j})$. So, $\sum\nolimits_{i=1}^n \sum_{j\in J} q_i\pi_{ij} = 1-\prod\nolimits_{j=1}^{n} (1-q_{j})$, which cannot exceed 1 since it is a probability.
\fin
\end{enumerate}

\noindent The following result shows that each customer is covered by its closest center. Its proof can be found in the Appendix.
\begin{teor}\label{closest}
There exists an optimal \ABC\ solution where every site is assigned to its closest center. Therefore, closest assignment constraints (CAC) can be used as valid inequalities.
\end{teor}
Observe that, in fact, the smaller assignment distances in a solution will seldom be the ones yielding the largest service cost.
Indeed, in order for this to happen, many other customers (those with larger assignment  distances) should have no demand. Therefore, the probability that a small assignment distance becomes the actual largest service cost can be extremely low.
For this reason, the approximation of the \ABC that only  accounts for the $K\le n$ largest assignment distances in the objective function can be very tight, even for moderate $K$ values (specially if probabilities $q_i$ are large).  From now on, we will refer to this approximation as \kABC. From a computational point of view, by using this approximation we avoid computing $\pi_{ij}$ probabilities associated with very small distances, that otherwise would require computing products of many demand probabilities, possibly causing stability and numerical problems. However, in contrast to the \ABC, now CAC are not automatically satisfied in general. Notice that they do hold in the homogeneous case, because in this case the resulting ordered median function has the isotonicity property \citep[see Section~3 and][]{nickelpuerto2005}.

\begin{lema}
In the K-PpCP, CAC must be explicitly included in the formulation. However they can be drop if all sites
share the same demand probability.
%In the \kABC, CAC must be explicitly included in the formulation, unless all sites share the same demand probability.
\label{lemaCAC-K}
\end{lema}

\begin{ejem}
\label{ejnoCAC}
Given the set of sites $N$ with coordinates $N=\{(81,65),(71,63),$ $(32,62),(22,72),$ $(70,21),$ $(44,34),(17,10),(25,36),(90,37),$
$(23,48)\}$, and using Euclidean distances, consider a \mbox{$3$-$P3CP$} instance with
\begin{figure}[!h]
\begin{center}
\vspace*{-6pt}
\fbox{\includegraphics[width=0.35\textwidth]{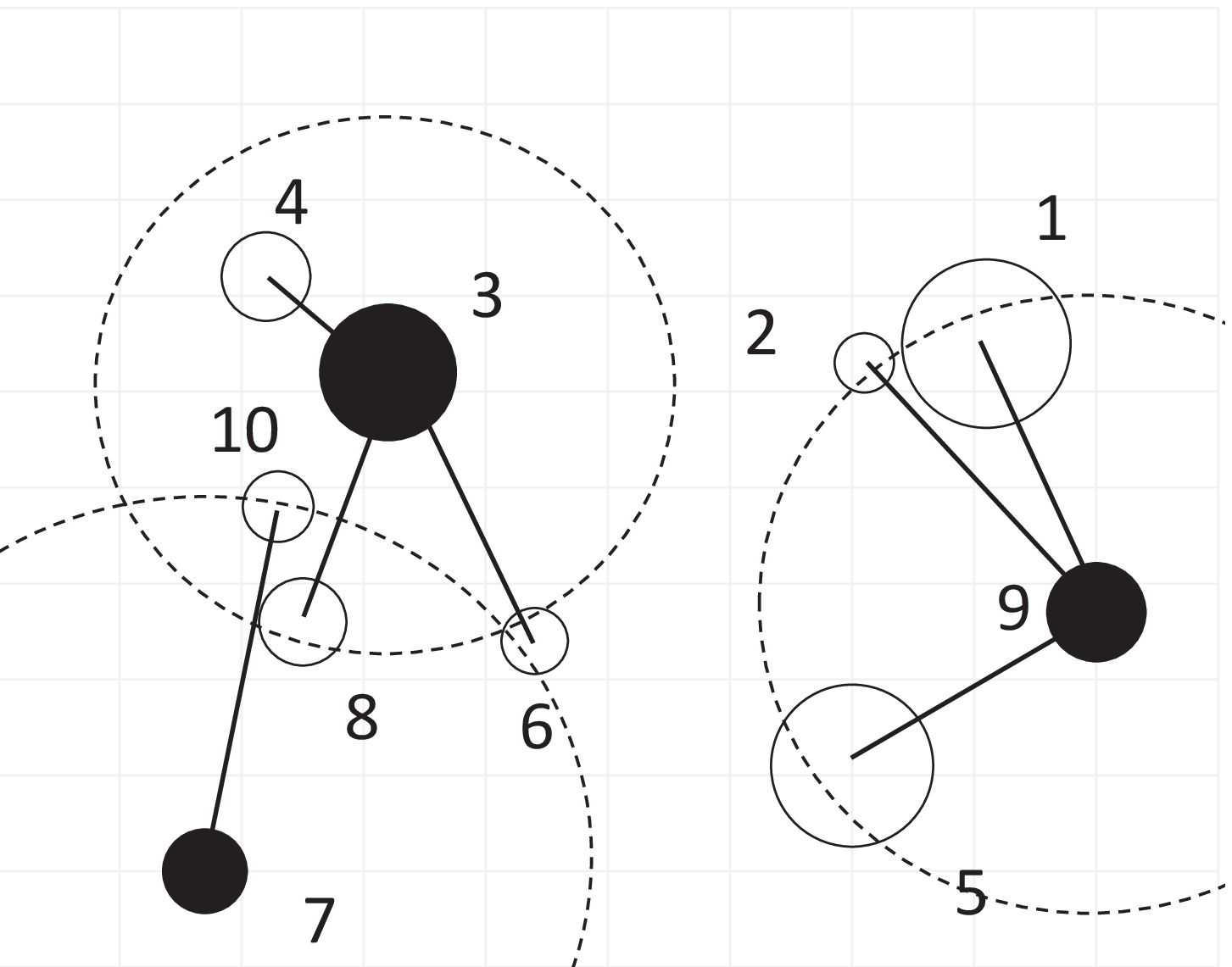}} \ \fbox{\includegraphics[width=0.35\textwidth]{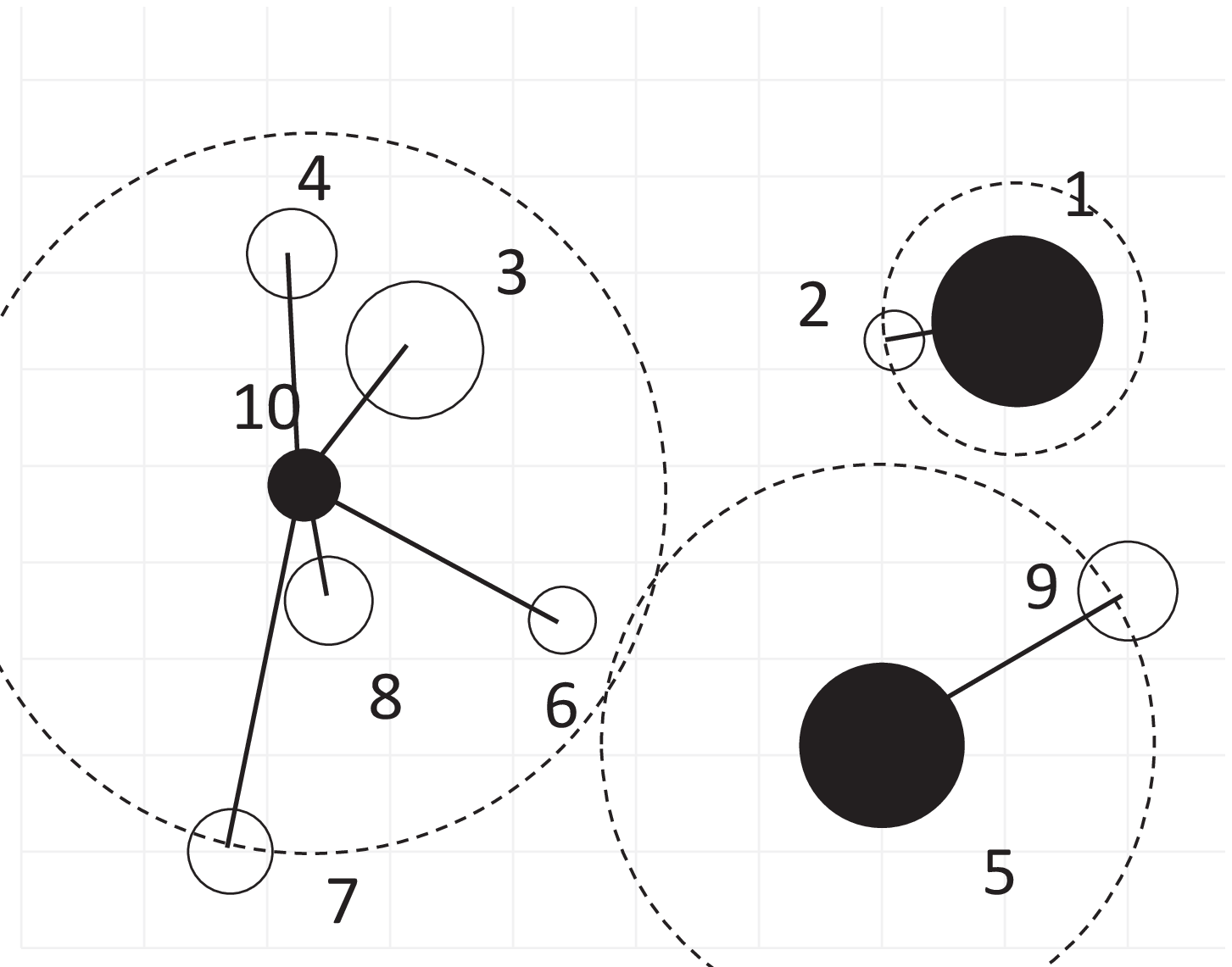}}
\end{center}
\caption{ Solutions of the instance in Example \ref{ejnoCAC} without CAC (left) and with CAC (right).}
\vspace*{-18pt}
\label{figure1}
\end{figure}
demand probabilities
$q=(0.97, 0.12, 0.63, 0.27,$ $0.9, 0.15, 0.24, 0.26, 0.33, 0.17)$. Figure~\ref{figure1} shows how, depending on whether CAC are imposed or not, the obtained solutions are different.
If CAC are not imposed, we obtain a solution with value $13.08$ (see Figure~\ref{figure1}, left). This solution allocates sites $4,\: 6$ and $8$ to center $3$, site $10$ to center $7$ and sites $1,\:2$ and $5$ to center $9$. However, in
this case the distance between site $10$ and the center located at $3$ ($d_{10, 3}=16.64$) is smaller than $d_{10,7}=38.47$.
By including CAC in the formulation, the objective value raises up to $17.58$ and the centers are located at sites $1$, $5$ and $10$. (See Figure~\ref{figure1}, right).
\end{ejem}

Note that the \ABC is equivalent to the \kABC with $K=n-p$ and, since $d_{ii}=0, \forall i\in N$ it makes no sense to take larger values of $K$. Therefore, in what follows we will present different formulations of the \kABC, for general $K\leqslant n-p$.

\section{Formulation for the homogenous case ($\mathbf{ q_i=q}$ for all $\mathbf{ i\in N}$)}
\label{section_homogeneous}
If all demand probabilities are equal, the probability that a given assignment provides the largest service cost depends on the number of larger assignment distances, but not on the associated customers. Therefore, the objective function of the homogeneous \kABC can be written as:
$$\sum_{k=n-K+1}^n q(1-q)^{n-k}d_{(k)}$$
where $d_{(k)}$ is the k-th value in the ordered assignment distances vector,
{\sl i.e.}, we are facing an ordered function. Thus, we can use the tools developed for discrete ordered median problems (DOMPs) \citep{nickelpuerto2005}.
 The formulation providing the best computational results for the DOMP is based on covering variables \citep{MNPV08}.
 However, the rationale behind this formulation cannot be adopted for the general \kABC. Indeed, unlike in the DOMP,
 additionally to the number of assignments with associated distances larger than a specific one, in the \kABC it is necessary to identify the customers defining those assignments. Since covering variables are based on the aggregation of equal assignment distances, they loose the  information on the customers defining them.
 For this reason, we next consider the three-index variables formulation, which can give better insights for the formulations we propose for the general \kABC that will be analyzed in the next sections.

Consider the set $T=\{n-K+1,\ldots,n\}$ and the binary variables:
\vspace*{-7pt}
\begin{itemize}
\itemsep1pt \parskip0pt \parsep0pt
\item For $i,j\in N, t\in T$, $x^t_{ij}$ takes value $1$ if and only if $i$ is allocated to $j$ and $d_{ij}$ is in the $t$-th position of the ordered assignment distances vector.
\item For $i,j\in N$, $x^{n-K}_{ij}$ is $1$ if and only if  $i$ is allocated to  $j$ and $d_{ij}$ is at position $t$, with $t\leqslant n-K$.
\end{itemize}
\vspace*{-7pt}
Additionally, we use the coefficients $\lambda^t = (1-q)^{n-t}$. The obtained formulation is:
\begin{alignat}{2}
\textrm{(FH)}\quad \min & \sum_{t=n-K+1}^{n} \lambda^t q \sum_{i=1}^n  \sum_{j=1}^n x^t_{ij} d_{ij}                      && \nonumber\\
s.t.                    &  \sum_{j=1}^n x_{jj}^{n-K}=p,                                                                  &&                    \label{cc_nfac.1} \\
	                    &\sum_{t=n-K}^n x_{ij}^t\leqslant x_{jj}^{n-K},                                                   &&\forall i,j\in N,   \label{cc_assign2b.1}\\
                        & \sum_{i=1}^n\sum_{j=1}^nx_{ij}^t=1,                                                             &&\forall t \in T, \label{cc_assign1.11}\\
                        &\sum_{i=1}^n\sum_{j=1}^n x_{ij}^{n-K}=n-K,                                                      &&                    \label{cc_assign1.1}\\
                        &\sum_{t=n-K}^n\sum_{j=1}^n x_{ij}^t=1,                                                           &&\forall i\in N,     \label{cc_assign2.1}\\
                        &\sum_{i=1}^{n}\sum_{j=1}^n d_{ij}x_{ij}^t\leqslant \sum_{i=1}^n\sum_{j=1}^n d_{ij}x_{ij}^{t+1},  &&\forall  t\in T\setminus\{n\},  \label{cc_sort.1}\\
                        &d_{ij}x_{ij}^{n-K}\leqslant \sum_{k=1}^n\sum_{l=1}^n d_{kl}x_{kl}^{n-K+1},                       &\quad &\forall i,j\in N,\label{cc_sort2.1}\\
                        & \sum_{t=n-K}^{n}{\sum_{\stackrel{a=1}{d_{ia}\succ  d_{ij}}}^n{x_{ia}^t}}+{x_{jj}^{n-K}}\leqslant 1,   &&\forall i,j\in N,                 \label{cac3ind2} \\
                        & x^t_{ij}\in\{0,1\} ,                                                                            && \forall i,j\in N,t\in T\cup\{n-K\}. \label{xtbinarias}
\end{alignat}
 Constraint \eqref{cc_nfac.1} guarantees that $p$ centers are located and constraints \eqref{cc_assign2b.1} ensure that each site is allocated to just one of them.
Constraints \eqref{cc_assign1.11} and \eqref{cc_assign1.1} guarantee that exactly one assignment takes the $t$-th position for $t\in T$ and the smallest $n-K$ assignment distances occupy the last positions. Constraints \eqref{cc_assign2.1} guarantee that each customer is associated with one position. The sorting of assignment distances is made through constraints \eqref{cc_sort.1} and \eqref{cc_sort2.1}. Finally, constraints \eqref{cac3ind2} are CAC \citep{cac}. As mentioned above, these constraints are valid but, actually, they are only necessary for the general case (Observe that here, in the homogeneous case, the objective function weights the ordered assignment distances with factors $q(1-q)^{(n-t)}$ that are monotonously increasing and, therefore, it has the isotonicity property). Note that, in case of ties among assignment distances of different customers, they can be sorted arbitrarily since all choices yield the same objective value. This formulation has served as a basis for the first formulations for the general case of the \kABC.

\section{Formulations for the ${\mathbf \kABC}$}
\label{sec_formulations}
\subsection{Three index formulation}
\label{sec_3ind}

In the general case $\lambda^t$ values are no longer known beforehand. Thus, they need to be replaced with decision variables. These new variables are defined as follows:
\begin{itemize}
\item For $i,j\in N,t\in T$, $\lambda_{ij}^t$ is the probability that there is not a service cost greater than $d_{ij}$, if $d_{ij}$ is in the $t$-th position of the ordered assignment distances, and $0$ otherwise. That is, $\lambda_{ij}^t = \pi_{ij}x_{ij}^t$.
\end{itemize}

 To force them take appropriate values, we will need some extra parameters. Let $q_{(1)}\leqslant \cdots \leqslant q_{(n)}$ be a nondecreasing sequence of the demand probabilities. We define $\kappa^t=\prod_{k=1}^{n-t} (1-q_{(k)})$. Then, the following formulation for the \kABC can be derived:
\begin{alignat}{2}
(\textrm{F3}^K)\quad \min \quad & \sum_{t=n-K+1}^{n}\sum_{i=1}^{n}\sum_{j=1}^{n} ( d_{ij}q_i)\lambda^t_{ij} &      & \label{cc_obj.1} \\
s.t.\quad & \textrm{constraints \eqref{cc_nfac.1}-\eqref{xtbinarias},}                                                     &&\nonumber\\
          & \lambda_{ij}^t \leqslant  \kappa^t x_{ij}^t ,                                                                   && \forall i,j\in N,t \in T, \label{cc_lambdax}\\
          & \sum_{i=1}^n \sum_{j=1}^n \lambda_{ij}^t = \sum_{i=1}^n \sum_{j=1}^n (1-q_i)\lambda_{ij}^{t+1},                 &\quad&\forall t\in
					T\setminus\{n\}, \label{lambdax2.1}\\
		  & \sum_{i=1}^n \sum_{j=1}^n \lambda^n_{ij} = 1,                                                                  && \label{lambdax3} \\
		  & \lambda_{ij}^t\geqslant 0,                                                                                      && \forall i,j\in N,t\in T. \label{lambdapos}
\end{alignat}

As explained in the last section, constraints \eqref{cc_nfac.1}-\eqref{cac3ind2} ensure that $x$
define properly sorted assignments. Now, as opposite to the homogeneous case, the sorting of equal-cost assignments can have an effect on the objective function value if ties occur between positions $n-K$ and $n-K+1$. In this case, we allow the least-cost ordering, which consists in assigning higher order to customers with lower demand probability. From now on, the order defined by $\prec$ will include this idea; {\sl i.e.}, if $i\neq i'$, $d_{ij}=d_{i'j'}$ and $q_{i'}> q_i$, we will consider that $d_{ij}\succ d_{i'j'}$. Constraints \eqref{cc_lambdax} ensure that $\lambda$ variables are consistent with the values of $x$ and  constraints \eqref{lambdax2.1}-\eqref{lambdax3} are used to compute the $\lambda$ variables.

\subsubsection{Valid inequalities}
\label{Valid_ineq_3ind}
\begin{itemize}
\item The probability that the largest service cost is among the $K$ largest assignment distances is
    \begin{equation}
    \sum_{i=1}^n\sum_{j=1}^n\sum_{t=n-K+1}^n q_i \lambda_{ij}^t\leqslant 1. \label{sumapi.1}
    \end{equation}
\item Combining \eqref{lambdax2.1} and \eqref{lambdax3} we obtain the valid equality:
    \begin{equation}
    \sum_{i=1}^n\sum_{j=1}^n\lambda_{ij}^{n-1} =\sum_{i=1}^n\sum_{j=1}^n x_{ij}^n(1-q_i). \label{uno}
    \end{equation}
\item The next inequalities are also valid:
    \begin{equation}
    \sum_{i'=1}^n\sum_{j=1}^n\lambda_{i'j}^{t-1} \geqslant\sum_{j=1}^n\sum_{ t'=n-K+1}^t(1-q_i)\lambda_{ij}^{t'}\quad \forall i\in N,t\in T\setminus\{n-K+1\}.
    \label{dos_prima}
    \end{equation}
        If $\sum_{j=1}^n\sum_{ t'=n-K+1}^t x_{ij}^{t'}=0$, the inequality holds trivially. Otherwise, if $x_{ij}^t=1$ for some $j\in N$ then, the corresponding \eqref{lambdax2.1} equation guarantees that \eqref{dos_prima} is satisfied. Again, due to \eqref{lambdax2.1}, if $x_{ij}^{t'}=1$ for some $n-K< t'<t$ then we have that $\lambda_{ij}^{t'}\leqslant\sum_{i'=1}^n\sum_{j'=1}^n\lambda_{i'j'}^{t-1}$ and \eqref{dos_prima} holds.

\end{itemize}

\subsubsection{Variable fixing }\label{3ind_fixing}

\noindent{\bf Trivial cases}

With the above assumptions, all sites where a center is located will be self-allocated, yielding the $p$ smallest assignment distances ($d_{ii}=0$). This allows
fixing to zero the $x_{ij}^t$ variables with:
\vspace*{-7pt}
\begin{itemize}
\itemsep1pt \parskip0pt \parsep0pt
\item $i\neq j$  and $ t\leqslant p$;
\item $i=j$  and  $t\geqslant \max\{p, n-K+1\}$; or
\item $t\geqslant n-K$ and $|\{j': d_{ij'} \prec d_{ij}\}|>n-p$.
\end{itemize}
\vspace*{-7pt}Clearly, the corresponding $\lambda_{ij}^t$ variables are automatically fixed to zero, too (by constraints \eqref{cc_lambdax}).

\noindent{\bf Fixing based on bounds}

The following lemmas provide some preprocesses that allow fixing some other $x$ and $\lambda$ variables.

\begin{lema}
\label{9.4}
Let $UB_{K-PpCP}$ be an upper bound of the \kABC. Then, if $ i,j\in N $ are such that  $q_i d_{ij}>UB_{K-PpCP}$, in any optimal solution $x_{ij}^t=0$, $\forall t\in T$.
\end{lema}
\dem
We will prove that for any feasible solution $X$ of the \kABC, with value $F_X$ we have that
\begin{equation*}
q_id_{ij}\leqslant F_{X}\quad \forall i,j\in N \textnormal{ \itshape{such that} }j\in X, d_{ij}=\min_{\ell\in X}\left\{d_{i\ell}\right\}.
\end{equation*}
Indeed, let $d_{i_nj_n}\geqslant\ldots\geqslant d_{i_1j_1}$ be the sorted list of assignment distances in $X$.

Then, $F_X=q_{i_n}d_{i_nj_n}+(1-q_{i_n})A_{n-1}$,
where $A_{n-1}= q_{n-1}d_{i_{n-1}j_{n-1}}+\sum_{s=n-K+1}^{n-2}{q_{i_s}\prod_{t=s+1}^{n-1}(1-q_{i_t})d_{i_sj_s}}$. Then, since $(1-q_{i_n})A_{n-1}\geqslant 0$, $F_X\geqslant q_{i_n}d_{i_nj_n}$. Moreover,
\begin{eqnarray*}
F_{X}=q_{i_n}d_{i_nj_n}+(1-q_{i_n})\left[q_{i_{n-1}}d_{i_{n-1}j_{n-1}}+(1-q_{i_{n-1}})A_{n-2}\right]&\geqslant& q_{i_n}d_{i_nj_n}+(1-q_{i_n})q_{i_{n-1}}d_{i_{n-1}j_{n-1}}\\
&\geqslant& q_{i_{n-1}}d_{i_{n-1}j_{n-1}}.
\end{eqnarray*}
The last inequality comes from the fact that $d_{i_nj_n}\geqslant q_{i_{n-1}}d_{i_{n-1}j_{n-1}}$. Accordingly, for $n-K< u\leqslant  n-2$,
\begin{eqnarray*}
F_{X} & \geqslant & q_{i_n}d_{i_nj_n}+\sum_{s=u+1}^{ n-1}{q_{i_s}d_{i_sj_s}\prod_{t=s+1}^n(1-q_{i_t})}+q_{i_{u}}d_{i_{u}j_{u}} \prod_{t=u+1}^n{(1-q_{i_t})}\\
&\geqslant& d_{i_{u+1}j_{u+1}}\left[ q_{i_n}+\sum_{s=u+1}^{ n-1} q_{i_s} \prod_{t=s+1}^n(1-q_{i_t})\right]+q_{i_u}d_{i_uj_u} \prod_{t=u+1}^n(1-q_{i_t})\\
&=&\left[1-\prod_{t=u+1}^n(1-q_{i_t})\right]d_{i_{u+1}j_{u+1}}+q_{i_u}d_{i_uj_u} \prod_{t=u+1}^n(1-q_{i_t}).\\
\end{eqnarray*}
and again, since $d_{i_{u+1}j_{u+1}}\geqslant q_{i_u}d_{i_u j_u}$, we have $F_X\geqslant q_{i_u}d_{i_u j_u}$.
Thus, taking $x_{\hat{\i} \hat{\j}}^t=1$  for a pair $\hat{\i}, \hat{\j} \in N$, such that $q_{\hat{\i}}d_{\hat{\i}\hat{\j}}>UB_{K-PpCP}$, and some $t\in T$, would yield a solution cost above $UB_{\kABC}$.
\fin
\begin{lema}
\label{9.5}
If $Ud^{t}$ is an upper bound on the $t$-th assignment distance,
$x_{ij}^{t'}=0\;\;\forall{i,j:d_{ij}>Ud^t; t'\leqslant t}$.
\end{lema}

\begin{lema}
\label{9.6}
If $Ld^{t}$ is a lower bound on the $t$-th assignment distance, $x_{ij}^{t'}=0\ \forall i,j:d_{ij}<Ld^t; t'\geqslant t$.
\end{lema}

\begin{lema}
\label{9.1}
The optimal value of the $pCP$ instance with distances $\tilde{d}_{ij}=q_{(1)}d_{ij}$ for  $i,j\in N$ (\mbox{$q_{(1)}=\displaystyle\min_{i\in N}{q_i}$}) yields a lower bound $\tilde{d}^*$ on the optimal \kABC value for any $K\geqslant 1$. Moreover, in any optimal solution,
$x_{ij}^n=0$ $\forall i,j\in N$ such that $d_{ij} q_i < \tilde{d}^*$.
\end{lema}
\dem
Let $X$ be an optimal solution of \kABC i.e., $X\subseteq\left\{1,\ldots,n\right\}$ and $|X|=p$. Using the notation of Lemma \ref{9.4}, its objective value is:
$F_X= q_{i_n}d_{i_nj_n}+\sum_{s=n-K+1}^{n-1}{q_{i_s}d_{i_sj_s}\left(\prod_{t=s+1}^{n}{(1-q_{i_t})}\right)}$.
 Hence, since $q_{(1)}\leqslant q_{i_n}$, we have that $q_{i_n}d_{i_nj_n}\geqslant q_{(1)}d_{i_nj_n}\geqslant \tilde{d}^*$.
\fin

\subsection{Compact 3-index formulation}
\label{sec_improved}

We next present a formulation that results from the aggregation of  variables used in the previous one. Together with the previous $\lambda$ variables, we  now consider:
\begin{equation}
\label{defast}
 x_{ij}=\sum\nolimits_{t=n-K}^n x_{ij}^t, \quad z_{it}=\sum\nolimits_{j=1}^n x_{ij}^t,\quad \forall i,j\in N\mbox{ and } t\in T,
\end{equation}
that allow building the following formulation:

\begin{alignat}{2}
\textrm{(CF3$^K$)}\  \min & \sum_{t=n-K+1}^{n}\sum_{i=1}^{n}\sum_{j=1}^{n} \lambda^t_{ij} q_id_{ij}&& \\
 s.t.\:
						  &\textrm{ constraints \eqref{lambdax2.1},}&&\nonumber\\
						& \sum_{j=1}^n x_{jj} \:=\: p                     ,                &&\label{p_abiertas} \\
						& x_{ij}\leqslant x_{jj,} &&\forall i,j \in N,\label{inst_abierta}\\
						& \sum_{i=1}^n  z_{it}   \:=\: 1,                  &&\forall t\in T, \label{constz2.1} \\
						& \sum_{j=1}^n  x_{ij}   \:=\: 1,       &&\forall i\in N, \label{asig_unica}\\
						& \sum_{\stackrel{a=1}{ d_{ia}\succ d_{ij}}}^n x_{ia}+x_{jj} \leqslant 1, &\qquad&\forall i,j\in N,\label{cac3indimp} \\
                        & \sum_{t=n-K+1}^n  z_{it}   \:\leqslant\: 1     ,       &&\forall i\in N\label{constz1.1} ,\\
                        & \lambda_{ij}^t\leqslant \kappa_t x_{ij}  , &&\forall i,j \in N,t \in T,\label{varrelation1}\\
                        & \sum_{j=1}^n\lambda_{ij}^t \leqslant z_{it}        ,      &&\forall i\in N, t\in T, \label{varrelation3}\\
                        & \sum_{k=1}^n \sum_{j=1}^n \lambda_{kj}^n d_{kj}\geqslant\sum_{j=1}^n x_{ij}d_{ij} ,&& \forall i\in N,\\
                        & (z_{it}+x_{ij}-1)t\leqslant\sum_{i_1=1}^n\sum_{\stackrel{j_1=1}{d_{i_1j_1} \preccurlyeq d_{ij}}}^n x_{i_1j_1},&\qquad&\forall i,j\in N,t \in T, \label{varrelation2}\\
						& x_{ij},z_{it} \:\in\: \{0,1\} ,               && \forall i,j\in N, t\in T, \\
                        & \lambda_{ij}^t \:\ge\: 0 ,                   && \forall t\in T.
\end{alignat}
Constraints \eqref{p_abiertas}-\eqref{constz2.1} are equivalent to \eqref{cc_nfac.1}-\eqref{cc_assign1.11}. Constraints
\eqref{asig_unica} and \eqref{constz1.1} ensure that each site is covered by only one center and takes one single position. CAC are given by \eqref{cac3indimp} where ties are treated as in $\textrm{F3}^K$ . Finally, constraints \eqref{varrelation1}-\eqref{varrelation2} ensure that $x$, $z$ and $\lambda$ take consistent values.

\begin{lema}
Integrality of assignment variables $x_{ij}$ with $i,j\in N, i\neq j$ can be relaxed.
\dem
If, for some $j\in N$ $x_{jj}=0$, then $x_{ij}=0$ for all $i\in N$ due to \eqref{inst_abierta}.
On the other hand,
if $x_{jj}=1$ and, for some $i,s \in N$, $x_{ss}=1$ and $d_{ij}\succ d_{is}$, by \eqref{cac3indimp}, we have that $x_{ij}=0$. Hence, by \eqref{asig_unica}, we have that $x_{ij}=1$ only if $x_{jj}=1$ and $x_{ss}=0\ \forall  s\in N: d_{is}\prec d_{ij}$.
\fin
\end{lema}
Now, the criteria presented in Section \ref{3ind_fixing} seldom allow to fix any $x$ variables. On the other hand, since CF3$^K$ uses the same $\lambda$ variables as before, they can be fixed using exactly the same criteria.

\subsubsection{Valid inequalities}
\begin{itemize} \itemsep1pt
\item If, in constraints \eqref{cc_sort.1} we replace $x_{ij}^t$ and $x_{ij}^{t+1}$ with $\lambda_{ij}^t$ and $\lambda_{ij}^{t+1}(1-q_i)$, we obtain:
\begin{equation}
\sum\limits_{i=1}^n \sum\limits_{j=1}^n d_{ij}\lambda_{ij}^t \leqslant \sum\limits_{i=1}^n \sum\limits_{j=1}^nd_{ij}(1-q_i)\lambda_{ij}^{t+1}\qquad \forall t \in T\setminus\{n\}. \label{impr_sort1_b} \end{equation}
\vspace*{-15pt}
\item Analogously to \eqref{uno}, using the definition of $z_{it}$ in \eqref{defast}, it holds that
\begin{equation}
\sum_{i=1}^n\sum_{j=1}^n\lambda_{ij}^{n-1}= \sum_{i=1}^n z_{in}(1-q_i).\label{uno_mejorada}
\end{equation}
\item Inequalities \eqref{lambdax3}, \eqref{sumapi.1} and \eqref{dos_prima} are also valid for this formulation.
\end{itemize}

\subsection{Formulation with probability chains}
\label{sec_probabchains}

In this section we adapt the formulation of the unreliable $p$-median problem proposed in \cite{Sergio} to the \kABC.
We denote $m=\frac{n^2+n}{2}$, the number of pairs $(i,j)$ such that $i,j\in N,\;i\leqslant j$  and
$M=\{1,\ldots,m\}$. Let $d'$ be the corresponding distances sorted in non-decreasing order (ties broken lexicographically). Also, we denote by $(i_k,j_k)$ the pair of sites associated with $d'_k$, $i_k\leqslant j_k$. Note that, for $k\leqslant n$, $i_k=j_k=k$ and $d'_k=0$. Now, we need the following variables defined for all $k\in M$.
\vspace*{-0.5em}
\begin{itemize}
\itemsep1pt \parskip0pt \parsep0pt
\item $y_k$ is the probability that the largest service cost is $d'_k$.
\item $\lambda_k$ is the probability that the largest service cost is $d'_{k'}$, with $k'<k$.
\item $s_k$, binary, takes value $1$ if and only if $d'_k$ is among the $n$-$K$ smallest assignment distances.
\end{itemize}
\vspace*{-0.5em}
We also use assignment variables $x_{ij}$ from formulation $CF3^K$. With all these variables we obtain:
\vspace*{-0.5em}
\begin{alignat}{2}
\textcolor{white}{.}\hspace*{-1em}\textrm{(PF$^K$)} \min &\; \sum_{k=1}^m d'_k y_k&& \\
								s.t.\: &\textrm{ constraints \eqref{p_abiertas}, \eqref{inst_abierta}, \eqref{asig_unica} and  \eqref{cac3indimp}},&&\nonumber \\
                                &y_m+\lambda_m=1,&&\label{rest3}\\
								&\lambda_k+y_k= \lambda_{k+1},&& \forall k\in M\!\!: k < m,\label{rest4}\\
								&y_k\leqslant q_{i_k}x_{i_kj_k}+q_{j_k}x_{j_ki_k},&& \forall k\in M\!\!: k>n,\label{rest5}\\
								&y_k\leqslant q_{i_k}x_{i_k i_k}, && \forall k\in M\!\!:k\leqslant n,\label{rest6} \\
								&y_k\geqslant q_{i_k}\lambda_{k+1}+x_{i_kj_k}-1-s_k, && \forall k\in M\!\!: k < m,\label{rest7}\\
								&y_k\geqslant q_{j_k}\lambda_{k+1}+x_{j_ki_k}-1-s_k, &&\forall k\in M\!\!: n < k < m,\label{rest8}\\
								&y_k\leqslant q_{i_k}\lambda_{k+1}+1- x_{i_kj_k}, &&\forall k\in M\!\!: k < m,\label{rest9}\\
								&y_k\leqslant q_{j_k}\lambda_{k+1}+1- x_{j_ki_k}, &&\forall k\in M\!\!: n < k < m,\label{rest10}\\
								&\sum_{k=1}^m s_{k}=n-K,&&\label{rest11}\\
								&s_k\leqslant x_{i_kj_k}+ x_{j_ki_k}, &&\forall k\in M\!\!: k>n,\label{rest12}\\
								& s_k=x_{i_ki_k}, &&\forall k\in M: k\leqslant n,\\
                                &Ks_k\leqslant \sum_{\stackrel{i,j=1}{ d_{ij}\succ d_{i_kj_k}}}^n  x_{ij} +K(1-x_{i_k j_k}), &&\forall k\in M,\label{rest14}\\
								&Ks_k\leqslant \sum_{\stackrel{i,j=1}{ d_{ij}\succ d_{j_ki_k}}}^n  x_{ij} +K(1-x_{j_k i_k}), &\qquad\qquad&\forall k\in M,\label{rest15}\\
								&s_k\in \{0,1\}, &&\forall k\in M,\\
								& x_{ij}\geqslant 0,\;x_{jj}\in  \{0,1\}, &&\forall i,j\in N.												
\end{alignat}

Constraints \eqref{rest3}-\eqref{rest10} guarantee the relationship
between $\lambda$ and $y$ variables to obtain consistent probabilities. Finally, constraints \eqref{rest11}-\eqref{rest15} ensure that $s$ variables
take the value $1$ only when the assignments associated with those variables are among the $n$-$K$ smallest distances.
Again, in case of ties between $d'_k$ and $d'_{k'}$, constraints \eqref{rest14} and \eqref{rest15}
consider that clients with smaller demand probabilities take higher positions.
 Notice that when there are no ties of a distance $d'_k$ with $k\in M$, \eqref{rest14} and \eqref{rest15} can be combined into the stronger constraint:
\begin{equation*}
Ks_k\leqslant \sum_{\ell>k}(x_{i_\ell j_\ell}+x_{j_\ell i_\ell}).
\end{equation*}

The following variables can be trivially fixed to zero:
\vspace*{-7pt}
\begin{itemize}
\addtolength{\itemsep}{-12pt}
\item $s_k$, if $k>m-K$,
\item $x_{ij}$, if $|\{j': d_{ij'} \succ d_{ij} \}|< p-1$.
\end{itemize}

Lemmas given in Section \ref{3ind_fixing} can also be adapted to this formulation to fix some of the $s$ variables.
In particular, Lemma \ref{9.5} can be applied to fix some of the $s$ variables to  $0$, by using an upper bound on the assignment distance occupying position $n-K$.
However, lemmas  \ref{9.4} and \ref{9.6}  now result on additional equations. For instance, if, for a given pair $(i_k,j_k)$ we could previously fix $x_{i_kj_k}^t$ to zero, for all $t\in T$, it means that either the assignment distance associated with $(i_k,j_k)$ is not incurred, or it is not among the $K$ largest ones. Therefore, in this case, this reasoning would not lead to fix to zero any variable in PF$^K$, but to set  $s_k = x_{i_kj_k}+x_{j_ki_k}$.

\section{Lower and upper bounds}\label{bounds1}
In this section we introduce some lower and upper bounds that will be used together with lemmas from
 Section \ref{3ind_fixing} to fix $x$ and $\lambda$ variables.

\begin{lema}\label{bound1}
The optimal solution of $pCP$ is an upper bound of the \kABC.
\end{lema}

\begin{teor}\label{teo1}
The solution of the following problem provides an upper bound for the \kABC.
\begin{alignat}{2}
\UB_1 = \min \quad & \sum_{t=n-K+1}^{n}\sum_{i=1}^{n}\sum_{j=1}^{n}  \kappa^{t} q_i d_{ij} x_{ij}^t   &      &\nonumber \\
s.t. \quad & \textrm{ constraints \eqref{cc_nfac.1}-\eqref{xtbinarias}}.  &&\nonumber
\end{alignat}
Recall that $\kappa^t=\prod_{k=1}^{n-t} (1-q_{(k)})$  (see Section \ref{sec_3ind}).
\end{teor}
\dem

Since $\kappa^t$ uses the
$n-t$ smallest probabilities, it bounds above the probability that none of the $n-t$ largest assignments is active. Consequently, $U\!B_1$ provides an upper bound on the \kABC.
\fin

%%%%%%%%%%%%%%%%%%%%%%%%%%%%TEOREMA 10.1%%%%%%%%%%%%%%%%%%%%%%%%%%%%%%%%%%%%%%%%%%
\begin{teor}\label{teo2}
The ordered median problem with weights $\lambda_t = q_{(n-t+1)}\kappa^t$ for $t\in T$, provides a lower bound for the \kABC. We will denote this bound with $\LB_1$.\\
{\rm The proof of Theorem \ref{teo2} is provided in the Appendix.}
\end{teor}
%%%%%%%%%%%%%%%%%%%%%%%%%%%%%%%%FINAL TEOREMA 10.1 %%%%%%%%%%%%%%%%%%%%%%%%%%%%%%%%%%%%%%%%%%%%%%%%%%%%%%%%

We denote the sorted sequence of distinct distances as
$0=d_{(1)}<\cdots< d_{(G)}=\displaystyle\max_{i,j\in N}\{d_{ij}\}.$
\begin{lema}\label{bound2}
For $h=1,\ldots, G$, consider the following problem:
\begin{alignat}{2}
n_U(h)= \max & \sum\limits_{(i,j): d_{ij}  \geqslant d_{(h)}} x_{ij} &&\\
                   s.t.  &\textrm{ constraints \eqref{p_abiertas}, \eqref{inst_abierta}, \eqref{asig_unica} and \eqref{cac3indimp}},&&\nonumber\\
									 &x_{ij}\in\{0,1\},&&i,j\in N. \label{xbin}
\end{alignat}
If $n_U(h)<n-t$, then $d_{(h)}$ is a strict upper bound on the $t$-th distance.
\end{lema}
\dem
$n_U(h)$ gives the maximum number of assignments that can be done at a distance not smaller than $d_{(h)}$ and is clearly non-increasing.
If $x_{ij}^t=1$ in a feasible solution, it means that $n-t$ assignments are made at distances $d_{ij}=d_{(h')}$ or larger,
so that $n_U(h')\geqslant n-t$ and  $n_U(h)\geqslant n-t$ for all $h\leqslant h'$.
\fin
\begin{lema}\label{bound3}
For $h=1,\ldots, G$, consider the following problem
\begin{alignat*}{2}
n_L(h)= \max & \sum\limits_{(i,j): d_{ij}  \leqslant d_{(h)}} x_{ij}, &&\\
                   s.t.  &\textrm{ constraints \eqref{p_abiertas}, \eqref{inst_abierta}, \eqref{asig_unica}, \eqref{cac3indimp} and \eqref{xbin}}.&&%\\
\end{alignat*}
If $n_L(h)<t-1$, then $d_{(h)}$ is a lower bound on the $t$-th distance.
\end{lema}
\dem
The same reasoning as before can be applied.
\fin
\begin{lema}\label{bound4}
Let $z_{p+t}$ be the optimal value of the $(p+t)CP$. Then $z_{p+t}$ is a lower bound on the $(n-t+1)$-th largest assignment distance of the \kABC. In particular, the optimal solution of the $(p+K)CP$ is a lower bound on any assignment distance.
\end{lema}
\dem
Let $X$ be the solution of the \kABC and $\left\{i_n,\ldots,i_{n-t+1}\right\}$ be the set of $t$ sites with the $t$-largest assignment distances . Then, $X\cup \{i_n,\ldots,i_{n-t+1}\}$ is a feasible solution of the $(p+t)CP$ with a cost that will not exceed $d_{n-t+1}$.
\fin

Finally, heuristic approaches can also be used in order to obtain upper bounds on the \kABC.
To this end, in Section \ref{sec_vns}, we adapt the VNS heuristic from \cite{dominguez2005}.

\section{Variable Neighborhood Search for the \kABC}
\label{sec_vns}
Variable Neighborhood Search (VNS) is a metaheuristic to solve combinatorial problems proposed by \cite{Mladenovic1997} for the $p$-median problem. It is a very well-known technique often used to solve discrete facility  location problems and it usually provides high quality solutions.
In particular, \cite{dominguez2005} and later \cite{pupega2014} proposed a VNS for solving the DOMP.% We will adapt this heuristic to solve the \kABC.

The VNS is based on a local search algorithm with neighborhood variations. Starting from a possible solution, the algorithm explores the neighborhoods in such a way that it obtains solutions progressively far from the current one.
In our problem, the $k$-th neighborhood is the set of feasible solutions that differ in $k$ centers from the current one.
Given a current solution, $x_{cur}$, characterized by a set of $p$ centers, $d_1(i)$ is the index of the center of $x_{cur}$ closest
to customer $i$ and  $d_{2}(i)$ is the index of the second closest center to customer $i$. Also, $f_{cur}$ is its objective value.

We use  an adaption of the algorithms described in \cite{dominguez2005} to our problem: Modified Move (MM), Modified Update (MU) and Modified Fast Interchange (MFI).
Given $x_{cur}$ and a new facility $j_{in}\in N\setminus x_{cur}$ to enter in the solution, MM finds the best facility $j_{out}\in x_{cur}$ to get out from the solution. Once we have $j_{in}$ and $j_{out}$, MU modifies vectors $d_{1}$ and $d_{2}$, i.e., this algorithm updates the value of the closest and second closest center for each customer according with the new set of facilities.
Finally, MFI uses MM and MU recursively to obtain the best modification of $x_{cur}$ in the current neighborhood.  It must be noticed that, the $k$-th neighborhood associated with $x_{cur}$ is defined as
$N_{k}(x_{cur)}=\{x: x\textrm{ is a set of $p$ centers with }|x_{cur}\setminus x|=k \}$.

In MM and MFI, the updates of the objective values $f_{cur}$ are necessary. The main difference between our heuristic and the one described in \cite{dominguez2005} resides in the evaluation of this objective function. In our case, given a set of $p$ candidate locations, we create a vector $d_{cur}$ with all the corresponding assignment distances ($d_{cur}(i)=d_{i\;d_{1}(i)}$). To evaluate the objective function we sort the indices vector $(1,\ldots,n)$ by non-increasing values of $d_{cur}$.  Using the indices and positions of the $K$ largest assignment distances we can obtain the function value for $x_{cur}$. A scheme of the VNS for the \kABC is the following:

\begin{description}\addtolength{\itemsep}{-10pt}
\item[Step 1] Initialize $x_{cur}$ with a random selection of $p$ locations. Compute $d_{1}$, $d_{2}$ and $f_{cur}$.
\item[Step 2]  We take $k=1$ and repeat the following steps until $k=p$:
\vspace*{-10pt}
\begin{itemize}\addtolength{\itemsep}{-3pt}
\item Repeat $k$ times:\\ Take a random center to be inserted in the current solution. Using MM,
obtain the best location to remove from $x_{cur}$ in turn. Use MU to update $x_{cur}$, $d_{1}$, $d_{2}$  and $f_{cur}$.
\item Apply MFI to find a better solution than $x_{cur}$ in $N_{k}(x_{cur})$.
 If necessary, update $x_{cur}$, $d_{1}$, $d_{2}$, $f_{cur}$ and take $k=1$.
 \end{itemize}
\end{description}

\section{Computational experience}
\label{sec_comput}
This section is devoted to the computational studies of the formulations and bounds that we described along the paper. After a brief description of the
instances used, we first evaluate the fixing preprocesses used and then we compare the three studied formulations. All of them were implemented  in the commercial solver Xpress 7.7 using the modeling language Mosel \footnote{ See \url{http://www.maths.ed.ac.uk/hall/Xpress/FICO_Docs/mosel/mosel_lang/dhtml/moselref.html}}.
All the runs were carried out on the same computer with an Intel(R) Core(TM) i7-4790K processor with 32 GB RAM. We remark that the cut generation of Xpress was disabled to compare the relative performance of formulations cleanly.

The instances used in this computational experience are
based on the p-median instances from ORLIB\footnote{Electronically available at \url{http://people.brunel.ac.uk/~mastjjb/jeb/orlib/files/}} (pmed1, pmed2, pmed3, pmed4
and pmed5). From each of them, we extracted several distance
submatrices with $n$ ranging in $\{6,10,13,15,$ $20,15,30\}$ and we considered $p\in\{3,5,7,10\}$. Besides, we took $K$ about the $20\%$ of $n$. Probability vectors $q$ were randomly generated, taking values between 0.01 and 1 rounded to $2$  decimals.

In what follows, we report aggregated results of the different experiments. Detailed results can be found in the supplementary material.

\subsection{Quality of the bounds} \label{sec:boundsQ}
We next evaluate the quality of the bounds on the \kABC presented in Section \ref{bounds1}.
Table \ref{tab11} shows, for instances of the same size, the average gap between each bound and the optimal solutions, and the CPU time (in seconds) required to compute them.
\begin{table}[!h]
\centering
\caption{Bounds: Average gaps and computing times\label{tab11}}
  \scalebox{0.9}{
\begin{tabular}{|rr|rr|rr|rr||rr|}
\hline
& & \multicolumn{2}{c|}{$\UB_1$} & \multicolumn{2}{c|}{VNS} & \multicolumn{2}{c||}{$pCP$} & \multicolumn{2}{c|}{$\LB_1$} \\
    n       & $\sharp$     & \% gap &  time &  \% gap  &  time &\% gap & time & \% gap  &  time \\
\hline
 6 & 5 & 5.24 & 0.02 & 0.00 & 0.00 & 112.49 & 0.01 & 55.85 & 0.04\\
10 & 10 & 13.91 & 0.07 & 0.00 & 0.00 & 140.53 & 0.02 & 40.63 & 0.16\\
13 & 15 & 22.07 & 0.25 & 1.38 & 0.00 & 63.66 & 0.03 & 46.52 & 0.39\\
15 & 15 & 19.36 & 0.43 & 0.66 & 0.01 & 58.59 & 0.05 & 46.21 & 0.73\\
20 & 15 & 33.08 & 2.88 & 0.01 & 0.02 & 38.38 & 0.14 & 46.16 & 4.05\\
25 & 15 & 46.06 & 34.10 & 0.66 & 0.03 & 31.60 & 0.33 & 45.07 & 15.98\\
30 & 15 & 52.88 & 246.75 & 1.04 & 0.07 & 18.60 & 0.68 & 48.49 & 99.92\\
\hline
\end{tabular}%
}
\end{table}%
The lower bound $\LB_1$ proved to be rather poor, with gaps close to $50\%$. Moreover, its computational burden increases very fast with the instance size. Regarding the upper bounds, it becomes evident that VNS provides the best results. Not only it yields the smallest gaps, which did not reach $1.5\%$ in any of the instance groups, but also the computational effort is very small (the whole set of instances was solved in less than 3 seconds in total).

Since VNS provides the best bounds with a small computational effort, we wanted to test it for larger instances. To this end, we generated a set of larger instances with $n\in \{50, 60, 70, 80\}$, $p=10$ and $q\in \{0.25, 0.5, 0.75\}$ from pmed7, pmed12, pmed17 and pmed22.  In order to be able to compare the obtained solutions with the optimal value, in this case we only considered homogeneous instances, which, as mentioned above, fit the structure of the DOMP. We implemented the formulation of the DOMP from \cite{MNPV08} and we run it with a time limit between 2 and 8 hours, depending on the instance size. The obtained results are given in Table \ref{tab10}.
\begin{table}[!h]
\caption{VNS for the homogeneous case}
  \centering
  \scalebox{0.9}{
    \begin{tabular}{|rrr|ccr|ccr|ccr|}
    \hline
\multicolumn{3}{|c|}{}      & \multicolumn{3}{c|}{q=0.25} & \multicolumn{3}{c|}{q=0.5} & \multicolumn{3}{c|}{q=0.75} \\
    \hline
     \multicolumn{1}{|c}{n} & \multicolumn{1}{c}{p} & \multicolumn{1}{c}{K} &  \multicolumn{1}{|c}{ gap$_\textrm{\tiny B\&B}$ } & \multicolumn{1}{c}{gap$_\textrm{\tiny VNS}$} & \multicolumn{1}{c}{time} & \multicolumn{1}{|c}{gap$_\textrm{\tiny B\&B}$ } & \multicolumn{1}{c}{gap$_\textrm{\tiny VNS}$} & \multicolumn{1}{c}{time} & \multicolumn{1}{|c}{ gap$_\textrm{\tiny B\&B}$ } & \multicolumn{1}{c}{gap$_\textrm{\tiny VNS}$} & \multicolumn{1}{c|}{time} \\
	\hline
50  & 10 &   11  &  0.01  &  1.99 &   1.88  &  0.01  &  1.91 &  2.15  &  0.01  &  6.97 &   2.21\\
60  & 10 &   13  &  0.72  &  2.53 &   2.79  &  0.15  &  0.08 &  3.95  &  0.01  &  3.23 &  5.93\\
70  & 10 &   15  &  5.81  &  0.00 &   8.26  &  0.53  &  0.71 &  8.54  &  1.42  &  3.58 &  7.00 \\
80  & 10 &   17  &  9.87  &  0.72 &  13.30  &  6.66  &  2.23 & 11.34  &  1.48  &  1.68 &  11.39 \\
    \hline
    \end{tabular}%
		}
  \label{tab10}%
\end{table}%
Columns under heading gap$_\textrm{\tiny B\&B}$ report the average, over the 5 instances of the same size, of the branch and bound \%gap at termination. Columns under gap$_\textrm{\tiny VNS}$ report the obtained \%gaps with respect to the optimal or the best known solution.  Finally, the average CPU requirements of the VNS are reported in the third column of each group. The quality of the solutions provided by the VNS, although being always good, seems to slightly deteriorate for larger $q$ values but it is not affected by the instance size.  As for the CPU times,
they increase quite smoothly with the instance size.

Recall that our interest on proposing bounds is their usefulness to fix variables according to the results of Section \ref{3ind_fixing}. Table \ref{tab:fixed} shows the minimum, average, and maximum value, of the percentage of  variables that could be fixed in formulation $F3^K$ for the above instances with $n\leqslant 30$.
\begin{table}[!h]
  \centering
	\caption{ Pertentage of fixed variables in $F3^K$\label{tab:fixed}}%
 \scalebox{0.9}{
    \begin{tabular}{rr|*{5}{r}|c|*{3}{r}}
          &       & L\ref{9.4}  & L\ref{9.5}  & L\ref{9.6}  & L\ref{9.6}$^*$ & L\ref{9.1}  & All   &  no  2, 4  & no 2  & no 4 \\
    \hline
          & min   & 27.8  & 0.0   & 1.1   & 5.3   & 0.5   & 46.7  & 45.4  & 46.7  & 45.4 \\
    $x$     & average & 48.4  & 1.3   & 7.6   & 8.9   & 1.8   & 61.4  & 61.0  & 61.2  & 61.1 \\
          & max   & 64.6  & 3.6   & 19.0  & 14.8  & 7.6   & 71.0  & 70.9  & 71.0  & 70.9 \\
  \hline
          & min   & 39.6  & 0.0   & 1.4   & 6.6   & 0.6   & 65.1  & 64.7  & 64.7  & 65.1 \\
    $\lambda$ & average & 59.2  & 1.6   & 9.3   & 11.2  & 2.4   & 75.0  & 74.9  & 74.9  & 75.0 \\
          & max   & 80.7  & 4.6   & 21.8  & 22.2  & 11.4  & 88.7  & 88.7  & 88.7  & 88.7 \\
    \end{tabular}%
}
\end{table}%
It must be pointed out that we used the results from Section \ref{bounds1} to obtain the necessary bounds.
In particular, since, as we have just seen, VNS provides the best upper bounds for our problem, we represent in the table the percentage of fixed variables with Lemma \ref{9.4} using VNS. Besides, in these results,  lemmas \ref{9.5} and \ref{9.6} use the bounds on the distances given by lemmas \ref{bound2} and \ref{bound3}, respectively. An alternative bound for using with Lemma \ref{9.6} is the one provided by Lemma
\ref{bound4}. In the table, we denote it by  L\ref{9.6}$^*$. The table also reports the percentage of fixed  variables given by the result
of Lemma \ref{9.1}. Finally, in the 3 last columns of the table we summarize the results of the best performing combinations, which exploit all results except Lemma~\ref{9.5} and/or Lemma~\ref{9.1}.

We can observe that the result with the largest impact is Lemma \ref{9.4}, which allows to fix between $27.8\%$ and  $64.6\%$ of the $x$ variables, and between $39.6\%$ and $80.7\%$ of the $\lambda$s.
Combining it with all the other lemmas, we can increase these ranges to $46.7\%$-$71.0\%$ and $65.1\%$-$88.7\%$, respectively. Almost the same figures are obtained by ignoring Lemma~\ref{9.5}, Lemma~\ref{9.1}, or both of them.

As mentioned before, in the case of formulation CF3$^K$, we can fix exactly the same $\lambda$ variables as for F3$^K$, but neither $x$ nor $z$ variables are fixed in this case. Despite this fact, formulation CF3$^K$ still remains smaller than F3$^K$ in general. Indeed, only in $2$ of the $90$ considered instances, the number of $x$ and $z$ variables in CF3$^K$ was larger than the number of non-fixed $x$ variables in F3$^K$. On the average, the number of $x$ and $z$ variables in CF3$^K$ was about $60\%$ of the number of non-fixed $x$ variables in F3$^K$ and this percentage tends to increase for large $p$ values, but to decrease for larger instances.

Finally, Lemmas from Section \ref{3ind_fixing} can be also adapted with the aim of fixing some of the $s$ variables of  formulation PF$^K$. Table \ref{tab:s_fixed} reports the percentage of fixed variables in this case.
\begin{table}[!h]
\caption{Percentage of included valid inequalities (L\ref{9.4}, L\ref{9.6}, L\ref{9.6}$^*$) and fixed $s$ variables (L\ref{9.5}) in $PF^K$\label{tab:s_fixed}}%
  \centering
 \scalebox{0.9}{
    \begin{tabular}{r|*{4}{r}}
                 & L\ref{9.4}  & L\ref{9.5}  & L\ref{9.6}  & L\ref{9.6}$^*$   \\
    \hline
             min     & 21.9 &  7.2 &   0.0 &  0.0    \\
          average & 43.2 &  31.0&   1.1 &  4.4    \\
             max     & 67.8 &  69.2&   4.8 &  9.1    \\
    \end{tabular}%
}
  \end{table}%
 Here, column L\ref{9.5} shows the percentage of $s$ variables that Lemma \ref{9.5} fixes to $0$. Besides, columns under headings L\ref{9.4}, L\ref{9.6} and L\ref{9.6}$^*$ report the percentage of $s$ variables for which we add the valid equalities $s_{k}=x_{i_kj_k}+x_{j_ki_k}$ using the mentioned lemmas.
Recall that, in this case, Lemma \ref{9.5} is applied after using Lemma \ref{bound2} to identify an upper bound on the assignment distance that occupies the $(n-K)$-th position (Ud$^{n-K}$).  As in Table \ref{tab:fixed}, we use VNS to obtain upper bounds for Lemma \ref{9.4}. Besides, to apply Lemma \ref{9.6} we use the lower bound Ld$^{n-K+1}$ provided by  Lemma \ref{bound3}. Again, column under heading L\ref{9.6}$^*$ reports the percentage of constraints that Lemma \ref{9.6} adds using the bounds from Lemma \ref{bound4}.
 In summary, Lemma \ref{9.4} is the one with the largest percentage of  included valid inequalities and now
Lemma \ref{9.5} allows to fix a significant percentage of variables too. Once more, the contribution of lemma  \ref{9.6}  with either bound is marginal.

\subsection{Evaluation of the formulations}

In this section we analyze the results of the alternative formulations that we described in the paper and we examine their different variants.

\subsubsection{Three index formulation}

To evaluate the impact of the different enhancements proposed for the three index formulation of Section~\ref{sec_3ind}, we have tested seven different variants, which are defined
by the combination of valid inequalities and fixing criteria used, and also by the type of approach used to add the inequalities (cut and branch - C\&{}B - or branch and cut - B\&{}C).
\begin{table}[!h]
\caption{ Variants of formulation F3$^K$ \label{F3k_variantes}}
 \centering
 \scalebox{0.8}{
 \begin{tabular}{rr|*{7}{ c}}
             &                   & F3$^K$:1   & F3$^K$:2                       & F3$^K$:3                       & F3$^K$:4                       & F3$^K$:5                       &  F3$^K$:6                       & F3$^K$:7\\
\hline
             & \eqref{sumapi.1}  & \checkmark & \checkmark                     & \checkmark                     & \checkmark                     & \checkmark                     & \checkmark                     & \checkmark                     \\
val. ineq.   & \eqref{uno}       &            &                                &                                &                                & \checkmark                     & \checkmark                     & \checkmark                     \\
             & \eqref{dos_prima} &            &                                &                                &                                &     C\&B                       &     B\&C                       &     C\&B                       \\
\hline
             & trivial           & \checkmark & \checkmark                     & \checkmark                     & \checkmark                     & \checkmark                     & \checkmark                     & \checkmark                     \\
             & L\ref{9.4}        &            & VNS                            & VNS                            & VNS                            & VNS                            & VNS                            & VNS                            \\
var. fixing  & L\ref{9.5}        &            & \checkmark                     &                                & \checkmark                     & \checkmark                     & \checkmark                     & \checkmark                     \\
             & L\ref{9.6}        &            & L\ref{bound3}, L\ref{bound4} & L\ref{bound3}, L\ref{bound4} & L\ref{bound3}, L\ref{bound4} & L\ref{bound3}, L\ref{bound4} & L\ref{bound3}, L\ref{bound4} & L\ref{bound3}, L\ref{bound4} \\
             & L\ref{9.1}        &            & \checkmark                     & \checkmark                     &                                & \checkmark                     & \checkmark                     &            \\
\end{tabular}%
}
\vspace*{-2ex}
\end{table}
Table~\ref{F3k_variantes} details the valid inequalities and criteria that have been considered in each variant. When both, C\&{}B and B\&{}C have been tested for the same family of valid
inequalities, the choice made in each variant is indicated in the table. In a similar way, the entry in the table indicates the bound used when different alternatives are available.
The decisions have been made according to the results of Section~\ref{sec:boundsQ} and preliminary computational tests.

The results for all these variants are summarized in Table~\ref{table_3ind1}. For each formulation variant, we report the LP gap (under ``Gap'') and the CPU time required to solve the instances (under ``Time''). Again, average values for equal sized instances are reported. In the cases where some of the $5$ instances remained unsolved after the time limit of $7200$ seconds, the number of such instances is provided in parenthesis next to the time and the average final gap is reported next to the LP gap. Also, the smallest time entry of each row is boldfaced.

\begin{table}[!h]
  \centering
		\caption{ Computational results for the three index formulation variants.  \label{table_3ind1}}
	\scalebox{0.68}{
    \begin{tabular}{|r@{/}c@{/}r|*{7}{r@{\;}r|}}
    \hline
 \multicolumn{3}{|c|}{}       & \multicolumn{2}{c|}{F3$^K$:1} & \multicolumn{2}{c|}{F3$^K$:2} & \multicolumn{2}{c|}{F3$^K$:3} &\multicolumn{2}{c|}{F3$^K$:4} & \multicolumn{2}{c|}{F3$^K$:5}
          & \multicolumn{2}{c|}{F3$^K$:6} &\multicolumn{2}{c|}{F3$^K$:7}\\
\hline
    n&p&K     & Gap    & Time  & Gap    & Time  & Gap    & Time  & Gap    & Time  & Gap    & Time  & Gap    & Time  & Gap    & Time \\
		\hline
    6     & 2     & 2     & 74.96 & \bf{0.06} & 60.42 & 0.15  & 60.42 & 0.07  & 61.61 & 0.16  & 60.39 & 0.06  & 60.42 & 0.15  & 61.59 & 0.09 \\
    10    & 3     & 3     & 79.47 & 0.78  & 34.43 & 0.78  & 34.43 & \bf{0.27} & 34.57 & 0.75  & 33.34 & 0.28  & 34.43 & 0.76  & 33.90 & 0.40 \\
    10    & 5     & 3     & 74.25 & 0.56  & 52.38 & 0.71  & 52.45 & 0.24  & 52.65 & 0.73  & 52.14 & \bf{0.24} & 52.38 & 0.73  & 52.43 & 0.43 \\
    13    & 3     & 4     & 89.94 & 8.66  & 50.61 & 2.72  & 50.63 & 1.30  & 50.62 & 2.67  & 50.21 & \bf{1.24} & 50.61 & 2.60  & 50.29 & 1.68 \\
    13    & 5     & 4     & 84.16 & 10.23 & 47.19 & 2.55  & 47.20 & \bf{1.17} & 47.19 & 2.57  & 45.72 & 1.37  & 47.19 & 2.61  & 45.73 & 1.62 \\
    13    & 8     & 4     & 65.06 & 1.50  & 40.55 & 1.69  & 40.55 & \bf{0.47} & 44.78 & 1.70  & 39.19 & 0.57  & 40.55 & 1.64  & 43.42 & 0.75 \\
    15    & 3     & 4     & 92.39 & 19.91 & 54.31 & 5.19  & 54.31 & \bf{2.38} & 54.32 & 5.19  & 53.78 & 2.70  & 54.31 & 5.30  & 53.80 & 3.05 \\
    15    & 7     & 4     & 83.57 & 25.72 & 47.91 & 4.37  & 47.95 & \bf{2.00} & 47.91 & 4.37  & 46.53 & 2.14  & 47.91 & 4.29  & 46.51 & 2.76 \\
    15    & 10    & 4     & 59.37 & 2.94  & 35.82 & 2.96  & 35.82 & \bf{0.74} & 35.82 & 2.96  & 35.06 & 0.85  & 35.82 & 3.05  & 35.30 & 1.26 \\
    20    & 3     & 5     & 95.19 & 410.76 & 53.91 & 26.27 & 53.91 & \bf{15.20} & 53.91 & 26.44 & 53.83 & 17.23 & 53.91 & 25.07 & 54.55 & 18.09 \\
    20    & 7     & 5     & 88.16 & 3100.67 & 40.57 & 50.99 & 40.60 & 31.47 & 40.57 & 52.19 & 40.21 & \bf{28.45} & 40.57 & 45.60 & 40.21 & 32.02\\
    20    & 10    & 5     & 82.86 & 962.15 & 41.41 & 23.93 & 41.45 & \bf{13.43} & 41.43 & 23.93 & 39.03 & 13.63 & 41.41 & 19.18 & 39.04 & 15.92 \\
    25    & 3     & 6     & 97.6{ $^{(21.6)}$} & 6390$^{(2)}$ & 48.55 & 179.50 & 48.56 & 124.32 & 48.55 & 180.67 & 48.42 & \bf{119.03} & 48.55 & 166.94 & 48.99 & 135.73 \\
    25    & 7     & 6     & { 95.4$^{(83.4)}$} & { 7201$^{(5)}$} & 42.53 & 680.36 & 42.53 & \bf{535.40} & 42.53 & 689.06 & 42.39 & 602.28 & 42.53 & 543.03 & 42.39 & 584.45 \\
    25    & 10    & 6     & { 93.0$^{(70.4)}$} & { 7200$^{(5)}$} & 42.10 & 688.08 & 42.10 & 587.91 & 42.10 & 699.31 & 41.48 & \bf{405.53} & 42.10 & 708.47 & 41.48 & 617.54 \\
    30    & 3     & 7     & { 98.6$^{(92.1)}$} & { 7204$^{(5)}$} & 49.91 & 1175.53 & 49.91 & 1247.29 & 49.91 & 1331.01 & 49.74 & \bf{1136.53} & 49.91 & 1435.94 & 50.20 & 1242.73 \\
    30    & 7     & 7     & { 95.9$^{(93.0)}$} & { 7203$^{(5)}$} & 52.9{$^{(18.2)}$} & 7294$^{(5)}$ & 54.0{ $^{(18.9)}$} & \bf{7101}$^{(4)}$ & 52.9{$^{(19.1)}$} & 7291$^{(5)}$ & 54.4{$^{(21.8)}$} & 7207$^{(5)}$ & 52.7{$^{(19.4)}$} & 7291$^{(5)}$ & 51.4{$^{(15.9)}$} & 7235$^{(5)}$ \\
    30    & 10    & 7     & { 89.1$^{(81.8)}$}	& { 7203$^{(5)}$} & 41.7{$^{(9.6)}$} & 5254$^{(3)}$ & 42.3{ $^{(12.0)}$} & 5299$^{(2)}$ & 42.1{$^{(10.4)}$} & 5419$^{(3)}$ & 40.1{ $^{(7.7)}$} & 5522$^{(2)}$ & 42.4{$^{(10.6)}$} & 5603$^{(3)}$ & 41.7{$^{(10.6)}$} & \bf{5050}$^{(3)}$ \\
    \hline
    \end{tabular}%
		}
\vspace*{-2ex}
\end{table}

Note that the number of $\lambda$ and $x$ variables fixed thanks to the results of Section~\ref{3ind_fixing} yield significant reductions of the computation times, allowing to increase the size of instances that can be solved. The variable fixing criterion provided by Lemma~\ref{9.5} does not improve on the combinations of the others. Indeed, among variants F3$^K$:2$-$4, the one excluding it (F3$^K$:3) seems to result in somehow smaller CPU times.
Note also that, since some variables are fixed according to optimality criteria, the LP gap is considerably reduced in variants F3$^K$:2$-$4 with respect to F3$^K$:1. However, they are still rather large and, unfortunately, the valid inequalities can only reduce them in some cases and by small amounts, resulting in similar times.

\subsubsection{Compact three index formulation}

As in the previous case, we have considered different alternative variants of formulation CF3$^K$, which are now detailed in Table~\ref{CF3k_variantes}. Trivial variable fixing has been applied in all cases.

\begin{table}[!h]
\caption{ Variants of formulation CF3$^K$ \label{CF3k_variantes}}
 \centering
 \scalebox{0.8}{
 \begin{tabular}{rr|*{5}{ c}}
             &                        & CF3$^K$:1  & CF3$^K$:2                     & CF3$^K$:3                    & CF3$^K$:4                   & CF3$^K$:5    \\
\hline
 &\eqref{lambdax3}, \eqref{sumapi.1}  & \checkmark & \checkmark                    &  \checkmark                 & \checkmark                   & \checkmark   \\
             & \eqref{dos_prima}      &            &                               & \checkmark                   &  C\&B                       &    C\&B      \\
val. ineq.   & \eqref{impr_sort1_b}   &            &                               &                              &                             & \checkmark    \\
             & \eqref{uno_mejorada}   &            &                               &  \checkmark                  & \checkmark                  & \checkmark    \\
\hline
             & L\ref{9.4}             &            & VNS                           &   VNS                        &     VNS                     &    VNS        \\
var. fixing  & L\ref{9.5}, L\ref{9.1} &            & \checkmark                    &  \checkmark                  &    \checkmark               & \checkmark    \\
             & L\ref{9.6}             &            & L\ref{bound3}, L\ref{bound4}  & L\ref{bound3}, L\ref{bound4} &L\ref{bound3}, L\ref{bound4} &         L\ref{bound3}, L\ref{bound4} \\
\end{tabular}%
}
\end{table}

Table \ref{table_3indimp} reports the results obtained with these variants of the compact three index formulation, on the same instances as before. The structure of
the table is exactly the same as for Table~\ref{table_3ind1}. Note that now, the lemmas of Section~\ref{3ind_fixing} are only used to fix $\lambda$ variables.
From Tables~\ref{table_3ind1} and \ref{table_3indimp} we can see that the LP bounds of the plain formulation CF3$^K$ are even looser than those of F3$^K$, although after applying all variable fixing criteria, the LP gaps become very similar in both formulations. Now, the inclusion of valid inequalities does have some mild impact on the CPU times required to solve the instances.

\begin{table}[!h]
  \centering
	   \caption{ Computational results for the compact three index formulation variants. \label{table_3indimp}}
  \scalebox{0.7}{
\begin{tabular}{|r@{/}c@{/}r|*{5}{rr|}}
    \hline
 \multicolumn{3}{|c|}{}       & \multicolumn{2}{c|}{CF3$^K$:1} & \multicolumn{2}{c|}{CF3$^K$:2} & \multicolumn{2}{c|}{CF3$^K$:3} &\multicolumn{2}{c|}{CF3$^K$:4} & \multicolumn{2}{c|}{CF3$^K$:5}\\
\hline
    n&p&K     & Gap    & Time  & Gap    & Time  & Gap    & Time  & Gap    & Time  & Gap    & Time   \\
		\hline
     6     & 2     & 2     & 87.32 & \bf{0.03} & 60.92 & 0.15  & 60.89 & 0.09  & 60.89 & 0.06  & 60.89 & 0.06 \\
    10    & 3     & 3     & 91.52 & 0.35  & 41.10 & 0.78  & 39.81 & \bf{0.28} & 39.81 & 0.30  & 39.81 & 0.32 \\
    10    & 5     & 3     & 95.18 & 0.34  & 55.53 & 0.78  & 54.14 & 0.24  & 54.14 & \bf{0.22} & 54.11 & 0.24 \\
    13    & 3     & 4     & 94.99 & 3.07  & 55.92 & 3.56  & 53.63 & 1.64  & 53.63 & \bf{1.61} & 53.63 & 1.97 \\
    13    & 5     & 4     & 95.54 & 4.51  & 51.14 & 2.89  & 47.93 & \bf{1.53} & 47.93 & 1.61  & 47.91 & 1.59 \\
    13    & 8     & 4     & 96.71 & 2.63  & 46.43 & 1.91  & 42.72 & 0.66  & 42.72 & 0.65  & 42.71 & \bf{0.62} \\
    15    & 3     & 4     & 97.25 & 7.57  & 60.40 & 8.32  & 57.78 & \bf{4.19} & 57.78 & 4.66  & 57.75 & 4.88 \\
    15    & 7     & 4     & 98.11 & 22.32 & 51.85 & 5.09  & 49.66 & 2.99  & 49.66 & \bf{2.74} & 49.47 & 3.30 \\
    15    & 10    & 4     & 98.05 & 6.96  & 40.88 & 3.48  & 38.26 & \bf{1.11} & 38.26 & 1.23  & 38.01 & 1.25 \\
    20    & 3     & 5     & 98.11 & 201.57 & 61.48 & 146.24 & 59.79 & 88.71 & 59.79 & \bf{81.89} & 59.76 & 105.09 \\
    20    & 7     & 5     & 98.10 & 1408.43 & 47.54 & 88.76 & 44.51 & \bf{82.99} & 44.51 & 95.16 & 44.14 & 84.85 \\
    20    & 10    & 5     & 98.51 & 1773.99 & 47.86 & 45.83 & 42.86 & \bf{27.15} & 42.86 & 31.29 & 42.63 & 35.77 \\
    25    & 3     & 6     & 98.15 & 3161.84 & 59.54 & 1644.53 & 58.00 & \bf{555.35} & 58.00 & 689.19 & 57.96 & 1213.56 \\
    25    & 7     & 6     & 98.3{$^{(74.8)}$} & 7201$^{(5)}$ & 49.65 & 2318.91 & 47.13 & \bf{957.09} & 47.13 & 1795.21 & 47.00 & 2353.50 \\
    25    & 10    & 6     &  { 98.6$^{(83.2)}$ }     &  {  7202$^{(5)}$}     & 48.71 & 2008.20 & 45.50 & 1114.48 & 45.50 & \bf{980.34} & 45.05 & 1458.29 \\
    30    & 3     & 7     &  {  98.6$^{(86.0)}$ }    &  {  7207$^{(5)}$}     & 63.5{$^{(19.3)}$} & 6749$^{(3)}$ & 60.17 & \bf{5428.08} & 61.0{$^{(9.4)}$} & 6627$^{(3)}$ & 59.0{ $^{(10.0)}$} & 6449$^{(3)}$ \\
    30    & 7     & 7     &  {  98.5$^{(95.6)}$ }    &   { 7207$^{(5)}$}     & 57.9{$^{(26.7)}$} & 7296$^{(5)}$ & 59.0{$^{(25.9)}$} & \bf{7209}$^{(5)}$
		& 59.3{$^{(27.6)}$} & 7211$^{(5)}$ & 57.8{$^{(27.9)}$} & 7212$^{(5)}$ \\
    30    & 10    & 7     &  {  98.7$^{(95.8)}$ }    &   { 7205$^{(5)}$}    & 53.2{$^{(22.5)}$} & 5987$^{(4)}$ & 50.2{$^{(18.9)}$} & \bf{5909}$^{(4)}$
		& 50.1{$^{(19.2)}$} & 5954$^{(4)}$& 47.9{ $^{(17.6)}$} & 6019$^{(4)}$ \\
		\hline
    \end{tabular}%
		}
\end{table}%	

However, although having a smaller number of variables, none of the variants of this formulation allows to solve to optimality all the instances that could be solved with some of the F3$^K$ variants.

\subsubsection{Formulation with probability chains}

The results of the last formulation proposed in this paper are reported in this section.
In this case, as shown in Table~\ref{tab:s_fixed}, the adaptation of the results from Section~\ref{3ind_fixing}
allows fixing a smaller fraction of the variables. Additionally, no valid inequalities were identified for PF$^K$.
For this reason, in this case we only considered two formulation variants;  PF$^K$:1, where only trivial variable fixing
is applied, and PF$^K$:2, where all the other criteria for fixing variables are considered.
Table \ref{results_chains} reports the corresponding results, following the same structure as in the previous sections.

\begin{table}[!h]
  \centering
			\caption{ Computational times for the formulation with probability chains.}
	\scalebox{0.8}{
    \begin{tabular}{|r@{/}c@{/}r|*{1}{rr|}|r@{/}c@{/}r|*{1}{rr|}}
    \hline
     n     & p     & K   & \multicolumn{1}{c}{PF$^K$:1} & \multicolumn{1}{c||}{PF$^K$:2} &
     n     & p     & K & \multicolumn{1}{c}{PF$^K$:1} & \multicolumn{1}{c|}{PF$^K$:2}
			\\
    \hline
    6     & 2     & 2       & \bf{0.02}      & 0.12 &         20    & 3     & 5    & 10.06             & \bf{5.32}  \\
    10    & 3     & 3       & \bf{0.14}      & 0.22 &         20    & 7     & 5    & 30.35             & \bf{ 15.92}  \\
    10    & 5     & 3       & \bf{0.15}      & 0.25 &         20    & 10    & 5    & 27.93             & \bf{ 12.94}  \\
    13    & 3     & 4       & \bf{0.51}      & 0.58 &         25    & 3     & 6    & 70.17             & \bf{ 67.33}  \\
    13    & 5     & 4       & \bf{0.68}      & 0.89 &         25    & 7     & 6    & 336.48            & \bf{ 157.87}  \\
    13    & 8     & 4       & \bf{0.32}      & 0.61 &         25    & 10    & 6    & 899.75            & \bf{ 360.01 } \\
    15    & 3     & 4       & \bf{1.09}      & 1.19 &         30    & 3     & 7    & \bf{1504.98}      & 2082$^{(1)\dag}$ \\
    15    & 7     & 4       & 1.88           & \bf{1.42}&     30    & 7     & 7    & 3036.39           & \bf{1702.35}  \\
    15    & 10    & 4       & \bf{0.49}      & 0.90 &         30    & 10    & 7    & 6845$^{(4)*}$   & \bf{ 4944}$^{(1)\bullet}$  \\
    \hline
    \end{tabular}%
		}	
		
	Average termination gaps: $^* 65.5\%$, $^\dag 19.6\%$ and $^\bullet 8.0\%$.
	 \label{results_chains}%
\end{table}%

 In Table \ref{results_chains} the LP gaps have not been included because the LP solution value was always $0$ and, consequently, the LP gaps were $100\%$.
In spite of this, we can compare {the efficiency} of the PF$^K$ formulations regarding the CPU times. If some of the instances in a group remained unsolved after two hours, Table \ref{results_chains} gives the average gap at termination. It is remarkable how this formulation improves on the CPU times of the previous ones. Besides, if all variable fixing
criteria are used CPU times are still further reduced. Note that, in this case, all instances but two were solved within two hours.

\subsubsection{Comparison of formulations}

\begin{figure}[!h]
\begin{center}
\vspace*{-6pt}
\includegraphics[width=0.9\textwidth]{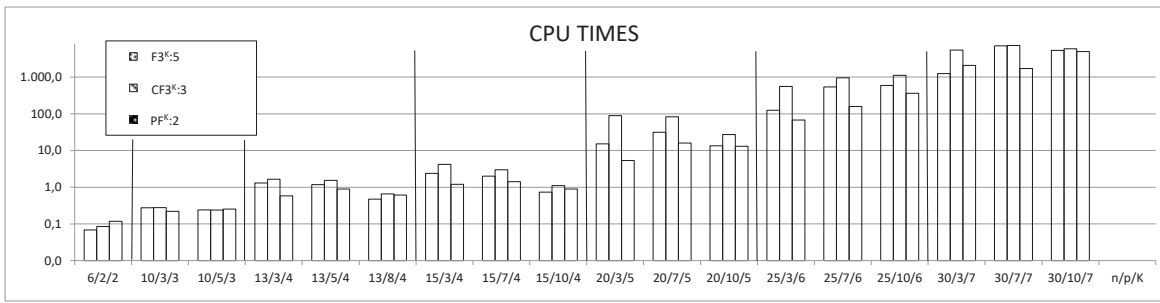}
\end{center}
\vspace*{-18pt}
\caption{ CPU times for the different variants\label{compare}}
\end{figure}

Following the results observed in the last subsections, we have chosen one representative variant of each formulation: F3$^K$:5 , CF3$^K$:3 and PF$^K$:2.
In order to compare them, Figure~\ref{compare} shows the times they yielded in logarithmic scale. Groups of instances with the same number of sites are delimited by vertical division lines.
The figure clearly shows that the computational burden of the \kABC grows exponentially with $n$ (mind the
 logarithmic scale in the vertical axis), but that this is specially true in the case of
formulation CF3$^K$. The superiority of PF$^K$ is evident here, although it can require the largest times
in some of the smaller instances. Moreover, recall that this formulation is the one that was able to solve the most
instances. Therefore the times of the unsolved instances for the other
two formulations are underestimated here and the figure shades the actual differences between them.
Within each group of instances with common $n$, we observe what usually happens with other
classical discrete location problems; they become more difficult as $p$ approaches $n/2$. The only possible
exception to this fact is formulation CF3$^K$, which tends to become more difficult as $p$ decreases.

 Summarizing, the PF$^K$ is the best formulation, since it allows to
solve most of the largest instances in the time limit. Besides, the adapted heuristic, VNS, provides accurate solutions for this generalization of the $p$-center problem in very small times.

\section{Concluding Remarks}
			
In this paper we introduce the probabilistic p-center problem (\ABC) and its generalization, the \kABC. These problems
allow to find compromise solutions, between the two extreme cases: the median-type problems, aimed at optimizing the average
service cost, and the center-type problems, aimed at optimizing the worst service level. To this end, the different sites
to serve are weighted according to their probability of requiring a service. In this way, one can prevent remote customers
with low demand probabilities from excessively conditioning the system configuration.

The particular case where all demand probabilities coincide fits in the structure of the ordered median problem and, therefore, it can be solved using all the tools available in the literature for it. However, for the general case, specific approaches need to be devised. The paper proposes and analyses three alternative formulations and a heuristic method.

Two of the formulations are based on existing formulations for the ordered median problem, while the third adapts some ideas that have been very successful for solving some reliable facility location problems. This last formulation dominates the other two. Given the superiority of this formulation based on probability chains, future research lines include the development of ad hoc procedures based on this formulation.

As for the proposed heuristic, it is an adaptation of a VNS procedure devised for the ordered median problem, and it provides high quality solutions in extremely reduced computation times.

\section*{Acknowledgements}
This research  has been partially supported by Ministerio de Econom\'ia y Competitividad under grants MTM2012-36163-C06-05,
MTM2013-46962-C02-02, MTM2015-63779-R, MTM2016-74983-C2-2-R and  by  FEDER-Junta de Andaluc\'ia under grant FQM 05849.

%\clearpage
%\newpage	
%\bibliographystyle{abbrvnat}
%\bibliographystyle{apalike}
%\bibliography{bibliografia}

\clearpage

\appendix
\section{Appendix}
\textbf{\underline{Proof of Theorem \ref{closest}}}:
Assume that two centers are located at sites $j$ and $j'$, and customer $i$ satisfies $d_{ij'}< d_{ij}$.
Consider solutions \textsc{Sol} where $i$ is covered by $j$, and \textsc{Sol'}, where $i$ is covered by $j'$, {\sl ceteris paribus}, and let $F_{j}$ and $F_{j'}$ denote their respective values. We will prove that  $F_j\geqslant F_{j'}$.

Let $d_{(1)}\leqslant \cdots \leqslant d_{(n)}$ and $d'_{(1)'}\leqslant \cdots \leqslant d'_{(n)'}$ be, respectively, the nondecreasing sequences of assignment distances in \textsc{Sol} and \textsc{Sol'}, and assume that $d_{ij}$ (resp. $d_{ij'}$) occupies position $t$ (resp. $s$) in its corresponding sequence.  By construction, $s\leqslant t$, and observe that $d'_{(s)'}= d_{ij'}$, $d_{(t)}=d_{ij}$,  $q'_{(s)'}=q_{(t)}=q_i$ and $d'_{(u)'}=d_{(u-1)}$  and $q'_{(u)'}=q_{(u-1)}$  for all $s+1 \leqslant u \leqslant t$.
\begin{eqnarray*}
F_{j'}-F_{j}&=&
\sum_{u=s}^{t}{\prod_{v=u+1}^n{(1-q'_{(v)})q'_{(u)}d'_{(u)}}}-\sum_{u=s}^{t}{\prod_{v=u+1}^{n}{(1-q_{(v)})}q_{(u)}d_{(u)}}\\
&=& \prod_{v=t+1}^{n} (1-q_{(v)}) \left[\sum_{u=s}^{t-1}q'_{(u)}d'_{(u)}\prod_{v=u+1}^t(1-q'_{(v)})+q'_{(t)}d'_{(t)}
-\sum_{u=s}^{t-1}q_{(u)}d_{(u)}\prod_{v=u+1}^t(1-q_{(v)})
-q_{(t)}d_{(t)}\right].
\end{eqnarray*}
To simplify the notation, let $F_{j'j}:=\frac{F_{j'}-F_{j}}{\prod_{v=t+1}^{n} (1-q_{(v)})}$. Then,
\begin{eqnarray*}
F_{j'j}&=&q'_{(s)}d'_{(s)}\hspace*{-1ex}\prod_{v=s+1}^t(1-q'_{(v)})+\hspace*{-1ex}\sum_{u=s+1}^{t-1}q'_{(u)}d'_{(u)}\hspace*{-1ex}\prod_{v=u+1}^{t}(1-q'_{(v)})+q'_{(t)}d'_{(t)}
-\hspace*{-1ex}\sum_{u=s}^{t-1}q_{(u)}d_{(u)}\hspace*{-1ex}\prod_{v=u+1}^{t}(1-q_{(v)})-q_{(t)}d_{(t)}\\
&=&q_{(t)}\left[d'_{(s)}\prod_{v=s}^{t-1}(1-q_{(v)})+\sum_{u=s}^{t-2}q_{(u)}d_{(u)}\prod_{v=u+1}^{t-1}(1-q_{(v)})+q_{(t-1)}d_{(t-1)}-d_{(t)}\right]\\
&\leqslant&q_{(t)}d_{(t)}\left[\prod_{v=s}^{t-1}(1-q_{(v)})+\sum_{u=s}^{t-2}q_{(u)}\prod_{v=u+1}^{t-1}(1-q_{(v)})+q_{(t-1)}-1\right] \leqslant 0.
\end{eqnarray*}
The last inequality is based on equation (\ref{eq:sumapis}).\fin

\noindent
\textbf{\underline{Proof of Theorem \ref{teo2}}:}
Let $X\subset N$ be the optimal solution of the \kABC and $F_X$ be its value.  Let $d_{n-K+1} \leqslant \cdots \leqslant d_n$ be the sorted list of the corresponding assignment distances involved in the objective function. To simplify the notation, and without loss of generality, we will assume that they correspond to sites $n-K+1, \ldots, n$, in this order.

\begin{eqnarray*}
F_X&=&q_nd_n+\sum_{t=n-K+1}^{n-1} q_t d_t \prod_{i=t+1}^n(1-q_i)\\
&=&q_nd_n+\ldots+\left[q_{n-i}d_{n-i}+(1-q_{n-i})q_{n-i-1}d_{n-i-1}\right] \prod_{s=n-i+1}^{n}(1-q_s)+\ldots+\\
&&+q_{n-K+1}d_{n-K+1}\prod_{t=n-K+2}^n(1-q_t).
\end{eqnarray*}

If there exists $q\in\{q_1,\ldots,q_n\}$ such that $q< q_{n-i}$ for $i< K$,
\begin{eqnarray*}
q_{n-i}d_{n-i}+(1-q_{n-i})q_{n-i-1}d_{n-i-1}&\geqslant& q d_{n-i}+(1-q)q_{n-i-1}d_{n-i-1},
\end{eqnarray*}
since
$q_{n-i}d_{n-i}+(1-q_{n-i})q_{n-i-1}d_{n-i-1}$ is an increasing function of $q_{n-i}$ and $q<q_{n-i}$.
Then,
\begin{eqnarray*}
F_X\hspace*{-0.8em}&\geqslant&\hspace*{-0.8em} q_nd_n+\ldots+\left[q d_{n-i}+(1-q)q_{n-i-1}d_{n-i-1}\right]\prod_{s=n-i+1}^{n}(1-q_s)
+\ldots
+q_{n-K+1}d_{n-K+1} \hspace*{-0.8em} \prod_{i=n-K+2}^n(1-q_i).
\end{eqnarray*}

This holds for all $i< K$. Consequently, if we define
\begin{eqnarray*}
F'_X&=&q^nd_n+\sum_{t=n-K+1}^{n-1}\prod_{i=t+1}^{n}(1-q^i)q^td_t,\textnormal{ with }q^{n-K+1},\ldots,q^n\in\{q_{(1)},\ldots,q_{(K)}\},
\end{eqnarray*}
\noindent where $q^i=q_i$ for any $i=n-K+1,\ldots,n$ if $q_i \in \{q_{(1)}, \ldots,q_{(K)}\}$, otherwise
$q^i$ is any element of $\{q_{(1)}, \ldots,q_{(K)}\}$, such that,
$\{q^{n-K+1}, \ldots,q^n\}=\{q_{(1)}, \ldots,q_{(K)}\}$. Then,
 we obtain that $F_X \geqslant F'_X$.
Since  $\{q^{n-K+1},\ldots,q^{n}\}=\{q_{(1)},\ldots,q_{(K)}\}$, there is a $q^{n-i}$ with $i\leqslant K$ such that $q^{n-i}=q_{(1)}$.
Then,
\begin{eqnarray*}
F'_X\hspace*{-0.8em}&=&\hspace*{-0.8em}q^nd_n+\ldots+\left[q^{n-i+1}d_{n-i+1}+(1-q^{n-i+1})q_{(1)}d_{n-i}\right] \hspace*{-0.8em}\prod_{s=n-i+2}^{n}\hspace*{-0.8em}(1-q^s)+\ldots
+q^{n-K+1}d_{n-K+1}\hspace*{-0.8em}\prod_{i=n-K+2}^n\hspace*{-0.8em}(1-q^i).
\end{eqnarray*}
We have that $q^{n-i+1}\geqslant q_{(1)}$ and $d_{n-i+1}\geqslant d_{n-i}$.  Then,
$d_{n-i+1}q^{n-i+1}+(1-q^{n-i+1})q_{(1)}d_{n-i}\geqslant d_{n-i+1}q_{(1)}+(1-q_{(1)})q^{n-i+1}d_{n-i}$.
As a result,
\begin{eqnarray*}
F'_X&\geqslant& q^nd_n+\ldots+\left[q_{(1)}d_{n-i+1}+(1-q_{(1)})q^{n-i+1}d_{n-i}\right]\prod_{s=n-i+2}^{n}(1-q^s)+\ldots+\\
&+&q^{n-K+1}d_{n-K+1}\prod_{i=n-K+2}^n(1-q^i).
\end{eqnarray*}
Following the same argument repeatedly,
$F'_X \geqslant q_{(1)}d_n+\ldots+(1-q_{(1)})q^{n-K+1}d_{n-K+1}\prod_{i=n-K+2}^{n-1}(1-q^i)$.
Since  $\{q^{n-K+1},\ldots,q^{n}\}=\{q_{(1)},\ldots,q_{(K)}\}$, there is a $q^{n-i}$ with $i\leqslant K$ such that $q^{n-i}=q_{(2)}$.
Then,
\begin{eqnarray*}
F'_X&\geqslant&q_{(1)}d_n+\ldots+(1-q_{(1)})\left[q^{n-i+1}d_{n-i+1}+(1-q^{n-i+1})q_{(2)}d_{n-i}\right]\prod_{s=n-i+2}^{n-1}(1-q^s)+\ldots+\\
&+&(1-q_{(1)})q^{n-K+1}d_{n-K+1}\prod_{i=n-K+2}^{n-1}(1-q^i).
\end{eqnarray*}
We have that $q^{n-i+1}\geqslant q_{(2)}$ and $d_{n-i+1}\geqslant d_{n-i}$. As before,
$d_{n-i+1}q^{n-i+1}+(1-q^{n-i+1})q_{(2)}d_{n-i}\geqslant d_{n-i+1}q_{(2)}+(1-q_{(2)})q^{n-i+1}d_{n-i}.$
Then, $F'_X\geqslant q_{(1)}d_n+(1-q_{(1)})q_{(2)}d_{n-1}+\ldots +(1-q_{(1)})(1-q_{(2)})q^{n-K+1}\prod_{i=n-K+2}^{n-2}(1-q^{i}).$
Following the same argument we can regroup $q_{(1)},\ldots,q_{(K)}$ and it holds
$$F_X\geqslant F'_X\geqslant q_{(1)}d_{n}+\ldots   +q_{(K)}d_{n-K+1}\prod_{i=1}^{K-1}(1-q_{(i)}). $$
\fin

\end{document}